\documentclass[pdflatex,sn-mathphys-num]{sn-jnl}%
\usepackage{silence}

\usepackage{pgfplots}
\usepackage{tikz-3dplot}
\usetikzlibrary{arrows.meta, calc, positioning}

\pgfplotsset{compat=1.18}

\usepackage{amsthm}

\newtheoremstyle{armatheorem}{12pt}{12pt}{\itshape}{0pt}{\bfseries}{.}{0.3em}{\thmname{#1}\thmnumber{ #2}\thmnote{ (#3)}}
\newtheoremstyle{armadefinition}{12pt}{12pt}{\normalfont}{0pt}{\bfseries}{.}{0.3em}{\thmname{#1}\thmnumber{ #2}\thmnote{ (#3)}}

\theoremstyle{armatheorem}
\newtheorem{armathm}{Theorem}[section]
\newtheorem{armalem}[armathm]{Lemma}
\newtheorem{armaprop}[armathm]{Proposition}

\theoremstyle{armadefinition}
\newtheorem{armadef}[armathm]{Definition}
\newtheorem{armarem}[armathm]{Remark}

\usepackage{mathptmx}
\usepackage[T1]{fontenc}
\usepackage{newtxtext}
\usepackage{anyfontsize}

\usepackage{titlesec}
\titleformat{\section}[block]
  {\centering\normalfont\normalsize\bfseries}
  {\thesection.}
  {0.5em}
  {}

\titleformat{\subsection}[block]
  {\centering\normalfont\normalsize\itshape}
  {\thesubsection.}
  {0.3em}
  {}

\usepackage{fancyhdr}
\usepackage{lastpage}
\usepackage{enumitem}

\usepackage{caption}

\DeclareCaptionFormat{armacaption}{%
    \centering\normalsize\textbf{#1}#2 \normalsize#3\par
}

\newcommand{\diff}{\mathop{}\!\mathrm{d}}
\newcommand{\beq}{\begin{equation}}
\newcommand{\eeq}{\end{equation}}

\newcommand{\ena}{\end{eqnarray}}
\newcommand{\bea}{\begin{eqnarray*}}
\newcommand{\eea}{\end{eqnarray*}}

\newcommand{\p}{{\partial}}

\newcommand{\h}{\hspace}

\newcommand{\interior}[1]{%
  {\kern0pt#1}^{\mathrm{o}}%
}

\newcommand{\normmm}[1]{{\left\vert\kern-0.25ex\left\vert\kern-0.25ex\left\vert #1
    \right\vert\kern-0.25ex\right\vert\kern-0.25ex\right\vert}}
\newtheorem*{theorem 2.1}{Theorem 2.1}

\def\XXint#1#2#3{{\setbox0=\hbox{$#1{#2#3}{\int}$}
  \vcenter{\hbox{$#2#3$}}\kern-.5\wd0}}

\usepackage{savesym}
\usepackage{amsmath}
\usepackage{amsfonts}

\savesymbol{Bbbk}

\usepackage[varg]{newtxmath}
\usepackage{newtxtext}

\usepackage{mathrsfs}

\usepackage{bm}

\makeatletter
\AtBeginDocument{
  \renewenvironment{proof}[1][\textit{Proof}]{\par
    \pushQED{\qed}%
    \normalfont\normalsize \topsep6\p@\@plus6\p@\relax
    \trivlist
    \item[\hskip\labelsep
          #1\@addpunct{.}]\ignorespaces
  }{%
    \popQED\endtrivlist\@endpefalse
  }
}
\makeatother

\usepackage{graphicx}%
\usepackage{multirow}%
\usepackage{amsmath,amssymb,amsfonts}%
\usepackage{amsthm}%
\usepackage{mathrsfs}%
\usepackage[title]{appendix}%
\usepackage{xcolor}%
\usepackage{textcomp}%
\usepackage{manyfoot}%
\usepackage{booktabs}%
\usepackage{algorithm}%
\usepackage{algorithmicx}%
\usepackage{algpseudocode}%
\usepackage{listings}%

\theoremstyle{thmstyleone}%

\theoremstyle{thmstyletwo}%

\theoremstyle{thmstylethree}%

\makeatletter
\renewcommand\normalsize{%
   \@setfontsize\normalsize{11bp}{13.2bp}%
   \abovedisplayskip 12\p@ \@plus2\p@ \@minus1\p@
   \abovedisplayshortskip \z@ \@plus3\p@%
   \belowdisplayshortskip 3\p@ \@plus3\p@ \@minus3\p@%
   \belowdisplayskip \abovedisplayskip%
   \let\@listi\@listI}
\renewcommand\small{%
   \@setfontsize\small{10bp}{12bp}%
   \abovedisplayskip 5\p@ \@plus3\p@ \@minus4\p@
   \abovedisplayshortskip \z@ \@plus2\p@
   \belowdisplayshortskip 3\p@ \@plus2\p@ \@minus2\p@
   \def\@listi{\leftmargin\leftmargini
               \topsep 4\p@ \@plus2\p@ \@minus2\p@
               \parsep 2\p@ \@plus\p@ \@minus\p@
               \itemsep \parsep}%
   \belowdisplayskip \abovedisplayskip}
\makeatother

\counterwithin{equation}{section}

\AtBeginDocument{\normalsize}
\begin{document}

\title[Biaxial Disclinations Around a Spherical Colloid]{\LARGE\itshape Biaxial Disclinations Around a Spherical Colloid}

\author{\fnm{Yuchen} \sur{Huang}\footnote{Email: 1155202952@link.cuhk.edu.hk, gentryhuang0403@gmail.com}\, \,
}

\author{\,  \fnm{Yong} \sur{Yu}\footnote{Email: yongyu@cuhk.edu.hk}}

\date{}
\affil{\orgdiv{Department of Mathematics}, \orgname{The Chinese University of Hong Kong}}

\abstract{We study the Landau–de Gennes theory in the one-constant limit for nematic liquid crystals in the exterior of a spherical colloid. The non-degenerate Nobili–Durand weak anchoring condition is imposed on the colloid boundary, supplemented by a homogeneous alignment condition at spatial infinity. Under an axially symmetric ansatz with an additional reflective symmetry and the Lyuksyutov constraint, we prove that biaxial boojums at the poles and a $-\frac{1}{2}$-degree Saturn ring disclination encircling the equator of the colloid can coexist within a single energy-minimizing configuration in this specific confined geometry. Finally, we rigorously characterize the local structural profile of the director field near these disclinations.}


\maketitle

\section{Introduction}\label{introduction}
Liquid crystals represent an intermediate phase of matter between conventional liquids and crystalline solids. The Landau–de Gennes theory \cite{de1993physics} provides a unified phenomenological framework for describing both uniaxial and biaxial nematic phases. Within this framework, the macroscopic orientational order is characterized by the \(Q\)-tensor—a real, symmetric, and traceless \(3 \times 3\) matrix. The local physical state of the liquid crystal is uniquely determined by the three real eigenvalues of \(Q\) and their relative magnitudes:\begin{itemize}
\item if all eigenvalues vanish, the phase is isotropic;\vspace{0.2pc}
\item if exactly two eigenvalues are equal and non-zero, the phase is uniaxial; \vspace{0.2pc}
\item if all three eigenvalues are distinct, the phase is biaxial.
\end{itemize}In a biaxial state, the preferred molecular alignment is dictated by the normalized eigenvector associated with the largest eigenvalue of \(Q\). This preferred orientation is conventionally referred to as the principal direction, or simply the director.

\subsection{Background}
The Landau–de Gennes model and its variants have been the subject of extensive mathematical and physical investigation. Crucially, classical liquid crystal theories describing the uniaxial director—such as the Oseen–Frank theory—can be recovered from the Landau–de Gennes framework through an appropriate limiting process \cite{majumdar2010landau}. Moreover, owing to the tensorial structure of its order parameter, the model can capture—both theoretically and numerically—a wide variety of topological defects observed experimentally. For instance, Alama, Bronsard, and Lamy \cite{alama2016minimizers} identified dipolar and quadrupolar disclinations within this framework in the limits of large and small colloidal particles, respectively. When a nematic droplet is immersed in an isotropic environment under confined geometries, three canonical types of topological defects emerge: the isotropic hedgehog point defect, the half-integer biaxial ring disclination, and the split-core segment disclination of strength 1. The physical structure of the radial hedgehog configuration has been extensively discussed \cite{gartland1999instability, kralj2001universal, rosso1996metastable, schopohl1988hedgehog}, and its stability has been mathematically established \cite{ignat2015stability, lamy2013some, majumdar2012radial}. For the half-integer ring disclination, early physical observations and numerical investigations were reported in \cite{kralj2001universal, lavrentovich1986phase, penzenstadler1989fine}, while the split-core segment disclination was first analyzed by Gartland and Mkaddem \cite{mkaddem2000fine}; the rigorous mathematical existence of both line defect configurations was later established by the second author and his collaborator \cite{tai2023pattern, yu2020disclinations}. These examples highlight the remarkable success of the standard Landau–de Gennes model in capturing both uniaxial and biaxial singularities in the director field. To accommodate even richer defect structures, Golovaty, Novack, and Sternberg \cite{golovaty2021novel} recently introduced a generalized Landau–de Gennes model featuring a quartic elastic energy density. For further examples and in-depth discussions on the Landau–de Gennes framework, we refer the reader to \cite{ambrosio1990existence, ambrosio1991boundary, an2017equilibrium, ball2010nematic, ball2011orientability, bauman2012analysis, bauman2016regularity, canevari2015biaxiality, canevari2017line, canevari2016radial, evans2016partial, henao2017uniaxial, ignat2016instability, majumdar2010equilibrium}.
\vspace{0.2pc}

In contrast, this article provides a mathematical analysis of biaxial surface point defects (boojums) and biaxial ring disclinations (Saturn rings) within the Landau–de Gennes framework. Colloids immersed in nematic liquid crystals inherently disrupt the underlying orientational order, generating topological defects in the director field \cite{liu2013nematic,muvsevivc2017liquid}. Under appropriate anchoring conditions, boojums and Saturn rings are commonly observed on the colloid surface and in the bulk, respectively. Boojums typically emerge at the antipodes of colloidal particles when molecules favor tangential surface alignment \cite{mermin1990boojums,mermin1977games,poulin1998inverted}; their core structures can be complex, exhibiting single-core, double-core, or split-core configurations \cite{tasinkevych2012liquid} (see \cite{muvsevivc2017liquid} for detailed schematics). Conversely, Saturn rings are disclination loops encircling the particle, driven by perpendicular alignment that induces a quadrupolar distortion. Under specific parameter regimes, these two defects can coexist—a phenomenon originally modeled in uniaxial frameworks by Volovik and Lavrentovich \cite{volovik1983topological, lavrentovich1998topological}, who demonstrated that such structures are induced by the topological dynamics of a bulk hedgehog defect as boundary anchoring transitions from perpendicular to tangential alignment (for representative numerical simulations, see \cite{FLTM2014}).\vspace{0.2pc}

Although numerous experimental observations and numerical simulations suggest the existence of boojums and Saturn rings and describe their formation mechanisms, rigorous mathematical verification remains elusive. Recently, Alama, Bronsard, and Golovaty \cite{ABG2019} made notable progress by identifying boojums in a two-dimensional thin-film model with oblique anchoring. However, the three-dimensional setting remains an open problem. This difficulty stems primarily from the complex interplay between boundary anchoring and bulk energy, as well as the multi-component structure of the $Q$-tensor order parameter, which complicates the rigorous analysis of boojum and Saturn ring core configurations. In this work, we employ the axially symmetric ansatz introduced by Yu \cite{yu2020disclinations} to reduce the model from a $Q$-tensor formulation to $\mathbb{S}^4$-valued axially symmetric harmonic maps (allowing for polynomial potential energy). By carefully analyzing the spatial relationship between the energy balls and the rotation axis, we establish Campanato decay estimates for the energy minimizers, thus proving the regularity of these $\mathbb{S}^4$-valued axially symmetric harmonic maps up to the colloid surface, except at finitely many singularities. We also provide mathematical justifications for the biaxial core structures of boojums and Saturn rings. Our work rigorously confirms the coexistence of boojum and Saturn ring disclinations reported by Dolganov et al. in \cite{dolganov2023statics}. We also refer the reader to a recent work by the second author and his collaborator \cite{taiyu2026}, where boojums are identified in a uniaxial model under a degenerate tangential boundary condition.

\subsection{Notations and Axial Symmetry}

Let $B_R \subset \mathbb{R}^3$ be a microscopic colloidal particle immersed in a liquid crystal medium, so that $G := \mathbb{R}^3 \setminus B_R$ represents the bulk domain occupied by the liquid crystal. Throughout this paper, the Cartesian, spherical, and cylindrical coordinates of $\mathbb{R}^3$ are denoted by $(x_1, x_2, x_3)$, $(r, \theta, \varphi)$, and $(\rho, \varphi, z)$, respectively, where $$r^2= x_{1}^{2}+x_{2}^{2}+x_{3}^{2},\h{15pt} \rho^2 = x_{1}^{2}+x_{2}^{2},\h{15pt} \theta\in \big[0,\pi\h{1pt}\big],\h{15pt} \varphi\in [-\pi,\pi],\h{15pt} z=x_{3}.$$ We denote the standard orthonormal basis in cylindrical coordinates by
\begin{align*}
\mathbf{e}_{\rho}\coloneq \left(\cos\varphi,\sin\varphi,0\right)^{\top},&\h{15pt} \mathbf{e}_{\varphi}\coloneq \left(-\sin\varphi,\cos\varphi,0\right)^{\top},\h{15pt} \mathbf{e}_{z}\coloneq  (0,0,1)^{\top}.
\end{align*}

The Landau–de Gennes energy functional equipped with a surface anchoring potential is given by
\begin{align*}
E_{a,W}(Q) \coloneq \int_{G} \Big( f_{E}(Q) + f_{B}(Q)-\min_{Q\in \mathcal{S}_{0}}f_{B} \Big) \, \mathrm{d}x + \int_{\partial B_R} f_{S}(Q)\, \mathrm{d}S.
\end{align*}
Here, $\mathcal{S}_{0}$ denotes the set of all real, symmetric, and traceless \(3 \times 3\) matrices, $f_E$ represents the elastic energy density under the one-constant approximation, $f_B$ is a fourth-order polynomial in $Q$ representing the bulk potential, and $f_S$ is a Rapini–Papoular type surface potential generalized to $Q$-tensors by Nobili and Durand \cite{muvsevivc2017liquid, nobili1992disorientation}. Explicitly, these density functions are given by
\begin{align*}
f_{E}(Q) &= \frac{1}{2} |\nabla Q|^2, \qquad f_{S}(Q) = \frac{W}{2} \big| Q - Q_{a}^{\circ} \big|^2, \\[5pt]
f_{B}(Q) &= -\frac{a^2}{2} |Q|^2 - \sqrt{6}\, \mathrm{tr}(Q^3) + \frac{1}{2} |Q|^4.
\end{align*}
The matrix norm on $Q$-tensors is induced by the Frobenius inner product, so that $\vert{}Q\vert{}^2 \coloneq \mathrm{tr}(Q^2)$. Furthermore, $-a^2$ denotes the reduced temperature, $W > 0$ is the anchoring strength, and $Q_a^\circ$ specifies the preferred target alignment on the boundary.\vspace{0.2pc}

Following \cite{tai2023pattern, yu2020disclinations}, we introduce an orthonormal basis for the linear space of $Q$-tensors:
\begin{equation*}
\begin{aligned}
     M_{1}=&\frac{1}{\sqrt{2}}\begin{pmatrix}
        0&0&1\\
        0&0&0\\
        1&0&0
    \end{pmatrix}, \h{15pt}  M_{2}=\frac{1}{\sqrt{2}}\begin{pmatrix}
        0&1&0\\
        1&0&0\\
        0&0&0
    \end{pmatrix},\h{15pt}
    M_{3}=\frac{1}{\sqrt{2}}\begin{pmatrix}
        0&0&0\\
        0&0&1\\
        0&1&0
    \end{pmatrix},\\[3mm]
    &M_{4}=\frac{1}{\sqrt{6}}\begin{pmatrix}
        -1&0&0\\
        0&-1&0\\
        0&0&2
    \end{pmatrix},\h{15pt} M_{5}=\frac{1}{\sqrt{2}}\begin{pmatrix}
        1&0&0\\
        0&-1&0\\
        0&0&0
    \end{pmatrix}.
\end{aligned}
\end{equation*}
Using this basis, we represent the order parameter $Q$ via the ansatz \begin{align}\label{axially ansatz of Q}
Q = Q[v]\coloneq \frac{a}{\sqrt{2}}\Big\{v_{1}(\cos 2\varphi \h{1pt} M_{5}+\sin 2\varphi \h{1pt} M_{2})+v_{2}M_{4}+v_{3}(\cos \varphi \h{1pt} M_{1}+\sin \varphi \h{1pt} M_{3})\Big\}.
\end{align}

\begin{armadef}[\bf Axial Symmetry]\leavevmode\begin{enumerate}\item[$\mathrm{(1)}$] A tensor field $Q$ given by \eqref{axially ansatz of Q} is called \textbf{axially symmetric} if $v_1, v_2, v_3$ depend only on $(\rho, z)$. \vspace{0.2pc}\item[$\mathrm{(2)}$] A $5$-vector field $w$ is called \textbf{axially symmetric} if it takes the form
\begin{align*}
    w = L[v] \coloneq (v_{1}\cos 2\varphi, \, v_{1}\sin 2\varphi, \, v_{2}, \, v_{3}\cos \varphi, \, v_{3}\sin \varphi)^{\top}
\end{align*}
for some $3$-vector field $v = v(\rho, z)$.\vspace{0.2pc}

\item[$\mathrm{(3)}$] An axially symmetric $5$-vector field $w = L[v]$ with $v = (v_1, v_2, v_3)^\top$ is called \textbf{$\mathcal{R}$-axially symmetric} (or simply \textbf{$\mathcal{R}$-symmetric}) if $v_1$ and $v_2$ are even with respect to $z$, whereas $v_3$ is odd with respect to $z$.
\end{enumerate}\end{armadef}

We now rescale the domain $G$ by setting $u(x) \coloneq v(Rx)$ for $x \in \Omega \coloneq \mathbb{R}^3 \setminus B_1$. Under this rescaling, $Q$ is an energy minimizer of $E_{a, W}$ under the axially symmetric ansatz if and only if $w = L[u]$ minimizes the rescaled functional $F_{a,\nu}$ under the axially symmetric ansatz, subject to suitable boundary conditions. Here, for a $5$-vector field $w$, the energy functional $F_{a,\nu}$ is defined by
\begin{align*}
F_{a,\nu}(w)\coloneq \int_{\Omega} \big|\h{.5pt} \nabla w \h{.5pt}\big|^{2}+\mu \left[D_{a}-3\sqrt{2}\h{1pt}S[w]+\frac{a}{2}(| w |^{2}-1)^{2}\right] +\nu\int_{\partial B_1} \big|\h{0.5pt} w - w_{a}^{\circ} \h{0.5pt}\big|^{2}.
\end{align*}
In this functional, $w_{a}^{\circ}$ denotes the smooth, axially symmetric $5$-vector field on $\partial B_1$ induced by $Q_{a}^{\circ}$. As $a \to \infty$, we assume that $w_{a}^{\circ}$ converges uniformly to a smooth unit $5$-vector field $w^{\circ}$ that inherits the same axial symmetry as $w_{a}^{\circ}$. In addition, we define $S[w]$ to be the $3$-homogeneous polynomial
\begin{align}\label{Sw polynomial}
S[w] \coloneq -w_{3}(w_{1}^{2} + w_{2}^{2}) + \sqrt{3} w_{2} w_{4} w_{5} + \frac{1}{2} w_{3}(w_{4}^{2} + w_{5}^{2}) + \frac{1}{3} w_{3}^{3} + \frac{\sqrt{3}}{2} w_{1}(w_{4}^{2} - w_{5}^{2}),
\end{align}
and introduce the non-dimensional constants
\begin{align*}
\mu \coloneq a R^2, \qquad \nu \coloneq W R, \qquad D_{a} \coloneq \frac{27}{16 a^3} \left[ 1 + \frac{4 a^2}{3} + \left( 1 + \frac{8 a^2}{9} \right)^{\frac{3}{2}} \right].
\end{align*}
Here, $D_a$ is a normalization constant chosen so that the minimum value of the potential density
\begin{align*}
f_{a}(w) \coloneq D_{a} - 3\sqrt{2} S[w] + \frac{a}{2} \big(|w|^2 - 1\big)^2
\end{align*}
is identically zero when restricted to the set $\big\{ L[v] : v \in \mathbb{R}^3 \big\}$. \vspace{0.2pc}

Near spatial infinity, we impose the far-field condition
\begin{equation}\label{far field f1}
w - H_a \h{0.5pt}e_3 \in \dot{H}^{1}(\Omega; \mathbb R^5),
\end{equation}
where
\begin{equation*}
e_3 \coloneq (0,0,1,0,0)^{\top} \quad \text{and} \quad H_{a} \coloneq \frac{3+\sqrt{9+8a^{2}}}{2\sqrt{2}a}.
\end{equation*}
$\dot{H}^{1}(\Omega; \mathbb R^5)$ denotes the homogeneous Sobolev space, defined as the closure of $C_c^\infty(\Omega \cup \partial B_1; \mathbb R^5)$ with respect to the norm
\begin{align*}
\| f \|_{\dot{H}^1(\Omega)} \coloneq \| \nabla f \|_{L^2(\Omega)}.
\end{align*}
Since $\dot{H}^1(\Omega; \mathbb R^5)$ is a Hilbert space under the inner product $\displaystyle \langle f, g \rangle_{\dot{H}^1} = \int_\Omega \nabla f : \nabla g$, it is weakly closed by Mazur's lemma.\vspace{0.2pc}

We define the admissible configuration space for $F_{a, \nu}$, denoted by $\mathcal{F}_{a, \nu}$, as
\begin{equation*}
\mathcal{F}_{a,\nu} \coloneq \bigl\{ w \in H_{\mathrm{loc}}^{1}(\Omega; \mathbb{R}^{5}) : w \text{ is axially symmetric, } F_{a,\nu}(w) < \infty, \text{ and } \eqref{far field f1} \text{ holds} \bigr\}.
\end{equation*}
Let $w_{a, \nu}$ be a minimizer of $F_{a, \nu}$ in $\mathcal{F}_{a, \nu}$. Up to a subsequence, it converges to some $w_\nu \in H_{\mathrm{loc}}^1(\Omega; \mathbb{S}^4)$, as $a \to \infty$. The convergence is strong in the space $H_{\mathrm{loc}}^1(\Omega; \mathbb R^5)$. Moreover, $\nabla w_{a,\nu}$ converges strongly to $\nabla w_{\nu}$ in $L^2(\Omega)$ as $a \to \infty$. Similarly to the finite-$a$ setting, we impose the following far-field condition for the limit map near spatial infinity:
\begin{align}\label{far field f2}
w - e_3 \in \dot{H}^{1}(\Omega).
\end{align}We can then show that $w_\nu$ is a minimizer of the limit energy
\begin{align}\label{first limit energy}
F_{\nu}(w) \coloneq \int_{\Omega}  |\nabla w|^2 + \sqrt{2}\mu \big(1 - 3 S[w]\big)  + \nu \int_{\partial B_1} |w - w^\circ|^2
\end{align}
within the configuration space
\begin{align*}
\mathcal{F}_{\nu} \coloneq \Big\{ w \in H_{\mathrm{loc}}^{1}(\Omega; \mathbb{S}^4) : w \text{ is axially symmetric}, \, F_{\nu}(w) < \infty, \text{ and } \eqref{far field f2} \text{ holds} \Big\}.
\end{align*}
Furthermore, it can be shown that $w_\nu$ satisfies the boundary value problem
\begin{equation}\label{EL of w}
\left\{\begin{aligned}
-\Delta w - \frac{3\sqrt{2}}{2} \mu \, \nabla_{w} S[w] &= \left( |\nabla w|^2 - \frac{9\sqrt{2}}{2} \mu S[w] \right) w && \text{in } \Omega, \\[2mm]
\partial_{n} w &= \nu \big[ w^\circ - (w^\circ \cdot w) w \big] && \text{on } \partial B_1.
\end{aligned}\right.
\end{equation}
Here, $\partial_{n}$ denotes the outward normal derivative along $\partial B_1$ with respect to the domain $\Omega$.\vspace{0.2pc}

Throughout this article, we write $(0,\infty)\times \mathbb{R}$ for the $(\rho,z)$-plane. Under this representation, the unit ball $B_1$ and its boundary $\partial B_1$ correspond to the domain $\mathbb{D}$ and its boundary $\partial \mathbb{D}$, respectively defined by\begin{align*}\mathbb{D} &\coloneq \left\{ (\rho, z) : z \in (-1, 1), \, 0 < \rho < \sqrt{1 - z^2} \h{1pt}\right\} \quad\text{and}\quad\partial \mathbb{D} \coloneq \left\{ (\rho, z) : z \in [-1, 1], \, \rho = \sqrt{1 - z^2} \h{1pt} \right\}.\end{align*}Let $D \coloneq (\partial_{\rho}, \partial_z)$ be the gradient operator on the $(\rho, z)$-plane. Assuming axial symmetry, the energy $F_{\nu}$ simplifies to\begin{align*}(2 \pi)^{-1} F_{\nu}(w) = E_{\nu}(u) \coloneq \int_{\overline{\mathbb{D}}^c} \left( |Du|^2 + \frac{1}{\rho^2} (4 u_1^2 + u_3^2) + \sqrt{2} \mu \big(1 - 3P(u)\big) \right) \rho \diff\rho \diff z + \nu \int_{-1}^1 |u - u^{\circ}|^2 \diff z,\end{align*}where $w^{\circ} = L[u^{\circ}]$ and $\overline{\mathbb{D}}^c$ denotes the complement of $\overline{\mathbb{D}}$ in $(0, \infty) \times \mathbb{R}$. The polynomial $P$ is defined by\begin{equation}\label{Polynomial P}P(v) \coloneq -v_{2} v_{1}^{2} + \frac{\sqrt{3}}{2} v_{1} v_{3}^{2} + \frac{1}{3} v_{2}^{3} + \frac{1}{2} v_{2} v_{3}^{2}.\end{equation}Note that when restricted to $\partial \mathbb{D}$, the map $u = u\big(\sqrt{1 - z^2}, z\big)$ depends solely on $z$. If $w_\nu = L[u_\nu]$ minimizes $F_\nu$ in $\mathcal{F}_\nu$, then $u_{\nu}$ minimizes $E_\nu$ in the configuration space $$\mathcal{E}_{\nu} \coloneq  \Big\{ u \in H^{1}_{\mathrm{loc}}\left(\overline{\mathbb{D}}^{\h{0.5pt}c}, \rho \diff\rho \diff z;  \h{1pt}\mathbb{S}^{2}\right) : E_{\nu}(u) < \infty \hspace{3pt}\text{and}\h{3pt} \eqref{far field f2}\,\,\,\text{holds for $L[u]$  } \Big\}$$ and satisfies the corresponding Euler–Lagrange equation\begin{align}\label{Euler-Lagrange equation of u} -\frac{1}{\rho}D \cdot  \big( \rho   D u \big) + \frac{1}{\rho^2} \begin{pmatrix}
   4u_{1}\\0\\u_{3}
\end{pmatrix} - \frac{3 \sqrt{2}}{2}\mu\nabla_u P(u)  =  \left\{ \big| D u\big|^2 + \frac{4 u_1^2 + u_3^2 }{\rho^2}  - \frac{9 \sqrt{2}}{2}\mu P\big(u \big) \right\} u.
\end{align}

\subsection{Main Results}

The primary goal of this paper is to construct a solution to \eqref{EL of w} that exhibits the coexistence of biaxial boojums and Saturn rings. The local structures of these disclinations are illustrated in Figure 1.  \tikzset{Boojum Ball/.style={x={(-0.353553cm, -0.353553cm)},y={(0.5cm, 0cm)},z={(0cm ,0.5cm)}}}
\begin{figure}[!htbp]
    \centering\scalebox{0.7}{
\begin{tikzpicture}[Boojum Ball]

\draw[->] (0,0,0)--(8,0,0) node [left]{$x$};
\draw[->] (0,0,0)--(0,8,0) node [right]{$y$};

\draw[-] (0,0,0)--(0,0,2.3);
\draw[dashed] (0,0,2.3)--(0,0,4);

\draw[->] (0,0,4)--(0,0,8) node [above]{$z$};

\draw[dashed] (0,0,0)--(0,0,-8);

\def\Rf{4.8}%
\def\Hf{1.1}%

\foreach \ang in {30,-30,60,-60,120,-120,-150,150,0,-180}{
    \pgfmathsetmacro{\u}{\Rf*cos(\ang)}
    \pgfmathsetmacro{\v}{\Rf*sin(\ang)}
    \pgfmathsetmacro{\xx}{-2.828427*\v}
    \pgfmathsetmacro{\yy}{2*(\u-\v)}
    \draw[gray!65,-{Stealth[length=2mm,width=1mm]}]
    (\xx,\yy,{-\Hf}) -- (\xx,\yy,\Hf);
}

\draw (0,-4,0) arc (180:360: 2cm and 0.707cm);
\draw[dashed] (0,-4,0) arc (180:0: 2cm and 0.707cm);
\draw (0,4,0) arc (0:50: 2cm and 2cm);
\draw (0,-4,0) arc (180:140: 2cm and 2cm);
\draw (0,-4,0) arc (180:360: 2cm and 2cm);
\draw[dashed] (0,-3.064,2.571) arc (140:50: 2cm and 2cm);
\node (1) at (0,4.2,0) [below] {$1$};

\draw [-{Stealth[length=2mm,width=1mm]}] (0,0.5,4)--(0,1.5,4);
\draw [-{Stealth[length=2mm,width=1mm]}] (0,-0.5,4)--(0,-1.5,4);
\draw [-{Stealth[length=2mm,width=1mm]}] (-0.5,0,4)--(-1.5,0,4);
\draw [-{Stealth[length=2mm,width=1mm]}] (0.5,0,4)--(1.5,0,4);
\draw [-{Stealth[length=2mm,width=1mm]}] (-0.5,-0.65,4)--(-1.3,-1.8,4);
\draw [-{Stealth[length=2mm,width=1mm]}] (0.5,0.65,4)--(1.3,1.8,4);
\draw [-] (2,-3.5,4)--(2,3.5,4)--(-2,3.5,4)--(-2,-3.5,4)--(2,-3.5,4) node [right] {};
\fill[black] (0,0,4) circle (2pt);
\node (2) at (0,0.35,4.7) {$\mathcal{N}$};

\draw [-{Stealth[length=2mm,width=1mm]}] (0,0.5,-4)--(0,1.5,-4);
\draw [-{Stealth[length=2mm,width=1mm]}] (0,-0.5,-4)--(0,-1.5,-4);
\draw [-{Stealth[length=2mm,width=1mm]}] (-0.5,0,-4)--(-1.5,0,-4);
\draw [-{Stealth[length=2mm,width=1mm]}] (0.5,0,-4)--(1.5,0,-4);
\draw [-{Stealth[length=2mm,width=1mm]}] (-0.5,-0.65,-4)--(-1.3,-1.8,-4);
\draw [-{Stealth[length=2mm,width=1mm]}] (0.5,0.65,-4)--(1.3,1.8,-4);
\fill[black] (0,0,-4) circle (2pt);
\node (2) at (-1,-0.7,-4) [] {$\mathcal{S}$};
\draw [dashed,-] (-2,1.7,-4)--(-2,-3.5,-4)--(-1.2,-3.5,-4);
\draw [-] (-1.2,-3.5,-4)--(2,-3.5,-4)--(2,3.5,-4)--(-2,3.5,-4)--(-2,1.7,-4) node [right] {};

\draw[black,dashed,line width=1pt] (0,6,0) arc (0:180: 3cm and 1.5cm);
\draw[black,line width=1pt] (0,-6,0) arc (180:360: 3cm and 1.5cm);
\draw[black,line width=1pt] (0,6,0) arc (0:59: 3cm and 1.5cm);
\draw[black,line width=1pt] (0,-6,0) arc (180:131: 3cm and 1.5cm);

\draw[] plot[domain=0:360, samples=60] (0, {6 + 1*cos(\x)}, {0 + 1*sin(\x)});
\fill[black] (0, 6, 0) circle (1pt);

\draw [blue,-{Stealth[length=2mm,width=1mm]}] (0,7,0)--(0,7,1);
\draw [blue,-{Stealth[length=2mm,width=1mm]}] (0,5.5,0.86)--(0,6.36,1.36);
\draw [blue,-{Stealth[length=2mm,width=1mm]}] (0,6.5,0.86)--(0,7,1.72);
\draw [blue,-{Stealth[length=2mm,width=1mm]}] (0,5,0)--(0,6,0);
\draw [blue,-{Stealth[length=2mm,width=1mm]}] (0,5.5,-0.86)--(0,6.36,-1.36);
\draw [blue,-{Stealth[length=2mm,width=1mm]}] (0,6.5,-0.86)--(0,7,-1.72);
\draw [blue,-{Stealth[length=2mm,width=1mm]}] (0,7,0)--(0,7,-1);

\draw[] plot[domain=0:360, samples=60] (0, {-6 + 1*cos(\x)}, {0 + 1*sin(\x)});
\fill[black] (0, -6, 0) circle (1pt);

\draw [blue,-{Stealth[length=2mm,width=1mm]}] (0,-7,0)--(0,-7,1);
\draw [blue,-{Stealth[length=2mm,width=1mm]}] (0,-5.5,0.86)--(0,-6.36,1.36);
\draw [blue,-{Stealth[length=2mm,width=1mm]}] (0,-6.5,0.86)--(0,-7,1.72);
\draw [blue,-{Stealth[length=2mm,width=1mm]}] (0,-5,0)--(0,-6,0);
\draw [blue,-{Stealth[length=2mm,width=1mm]}] (0,-5.5,-0.86)--(0,-6.36,-1.36);
\draw [blue,-{Stealth[length=2mm,width=1mm]}] (0,-6.5,-0.86)--(0,-7,-1.72);
\draw [blue,-{Stealth[length=2mm,width=1mm]}] (0,-7,0)--(0,-7,-1);

\fill[black] (0,0,7.2) circle (2.3pt);
\fill[black] (0,0,-6.4) circle (2.3pt);

\end{tikzpicture}}\vspace{0.2pc}
    \caption{Disclinations around a colloid}
    \label{Disclinations around the colloid}
\end{figure}
\noindent Specifically, two boojums are located at the north and south poles ($\mathcal{N}$ and $\mathcal{S}$), while the thick ring encircling the equator represents the Saturn ring disclination. The arrows indicate the corresponding director field. \vspace{0.2pc}

To present our main results, we first introduce some key concepts.
\begin{armadef}\label{vector field singularity}
A point $x \in \Omega$ is called a regular point of a vector field $v$ if $v$ is Lipschitz in a neighborhood of $x$. Otherwise, $x$ is called a singular point (or singularity) of $v$.\end{armadef}
\begin{armadef}
Let $u = u(\rho, z)$ be a unit $3$-vector field in $\Omega$. The following three quantities
\begin{equation}\label{eigenvalues} \begin{aligned}
\lambda_{1}\coloneq  -\frac{1}{2} \left( u_1 + \frac{1}{\sqrt{3}} u_2 \right), \h{15pt}
&\lambda_{2}\coloneq \frac{1}{4} \left\{\h{1pt} \left( u_1 + \frac{1}{\sqrt{3}} u_2 \right) - \sqrt{\big( u_1 - \sqrt{3} u_2 \big)^2 + 4 u_3^2 } \h{1pt} \right\},\\[2mm]
&\lambda_{3}\coloneq \frac{1}{4} \left\{ \h{1pt} \left( u_1 + \dfrac{1}{\sqrt{3}} u_2 \right) + \sqrt{\big( u_1 - \sqrt{3} u_2 \big)^2 + 4 u_3^2 } \h{1pt}\right\}
\end{aligned}
\end{equation} are called eigenvalues induced by $u$. They are eigenvalues of the matrix \begin{align}\label{matrix}
\frac{1}{\sqrt{2}}\Big\{u_{1}(\cos 2\varphi \h{1pt} M_{5}+\sin 2\varphi \h{1pt} M_{2})+u_{2}M_{4}+u_{3}(\cos \varphi \h{1pt} M_{1}+\sin \varphi \h{1pt} M_{3})\h{1pt}\Big\}.\end{align} At a point $x \in \Omega$, we say \begin{enumerate}[label=\textup{(\arabic*).}, leftmargin=*]
\item $u$ is isotropic if $x$ is a point singularity of $u$ in the sense given in Definition \ref{vector field singularity};\vspace{0.2pc}
\item $u$ is uniaxial if two of  the three eigenvalues in \eqref{eigenvalues} are identical and different from 0 at $x$; \vspace{0.2pc}
\item $u$ is biaxial if all the three eigenvalues in \eqref{eigenvalues} are different at $x$.\end{enumerate}
\end{armadef}

Using the three eigenvalues introduced in Definition \ref{eigenvalues}, we define the director field as follows:
\begin{armadef}
A vector field is called a director field determined by $u$ if it is a normalized eigenvector corresponding to the largest eigenvalue among $\lambda_1, \lambda_2$, and $\lambda_3$. Furthermore, if the director field can be defined smoothly in a neighborhood of $x$, it is said to be regular at $x$. Otherwise, $x$ is called a disclination (or a singularity) of the director field.
\end{armadef}

We also impose certain assumptions on the preferred alignment $w^\circ = L[u^\circ]$, which are stated below in terms of $u^\circ$: \begin{align}\label{assump 1}
u^{\circ} \equiv - (0,1,0)^{\top}\,\,\text{close to $\mathcal N$ and $\mathcal S$}, \qquad \big(u^{\circ}\big)_{1} \h{1pt} \big\vert_{\partial B_{1}}\geqslant 0, \qquad \big(u^{\circ}\big)_{1}\,\, \text{is not constant on $\p B_1$}.\end{align}

With these preparations, we now state our first result concerning the existence of boojums at $\mathcal{N}$ and $\mathcal{S}$.
\begin{armathm}\label{main theorem 1}
Suppose $w_{\nu} = L[u_\nu]$ minimizes the energy $F_{\nu}$ over $\mathcal{F}_{\nu}$. Then the following hold:
\begin{enumerate}
    \item [$\mathrm{(1).}$] For any $\nu > 0$, the minimizer $w_{\nu}$ is regular on $\overline{\Omega}$, except possibly at a finite number of singular points located on the $z$-axis, $l_z$. Furthermore, $w_{\nu}$ is regular in a neighborhood of each pole $\mathcal{N}$ and $\mathcal{S}$.\vspace{0.2pc}
    \item [$\mathrm{(2).}$] There exists $\nu_0 > 0$ such that $w_{\nu}(\mathcal{N}) = w_{\nu}(\mathcal{S}) = - e_3$ for each $\nu > \nu_0$. \vspace{0.2pc}
    \item [$\mathrm{(3).}$] For each $\nu > \nu_0$, there exists a sufficiently small $\epsilon_{\nu} > 0$ such that$$\lambda_{3} > \lambda_1 > \lambda_2 \quad\text{on }\left(B_{\epsilon_{\nu}}(\mathcal{N}) \cap \Omega \right)\setminus l_{z}.$$The same ordering holds in a neighborhood of $\mathcal{S}$.\vspace{0.2pc}
    \item [$\mathrm{(4).}$] In a neighborhood of each pole, $\mathcal{N}$ and $\mathcal{S}$, the director field is given by the normalized eigenvector corresponding to the leading eigenvalue $\lambda_3$. More precisely, it can be expressed as follows:
\begin{equation}\label{n-field formula}
\begin{split}
\mathbf{d}\left[u_{\nu}\right] \coloneq  &\frac{\sqrt{2}}{2} \left( 1 + \frac{u_{\nu, 1} - \sqrt{3}\h{0.5pt}u_{\nu, 2}}{\sqrt{\left(u_{\nu, 1} - \sqrt{3}\h{0.5pt}u_{\nu, 2}\right)^2 + 4u_{\nu, 3}^2}} \right)^{1/2} \mathbf{e}_{\rho} \\[2mm] &\h{10pt}+ \frac{\sqrt{2} \h{0.5pt}u_{\nu, 3}}{\sqrt{\left(u_{\nu, 1} - \sqrt{3}\h{0.5pt}u_{\nu, 2}\right)^2 + 4u_{\nu, 3}^2}}\left( 1 + \frac{u_{\nu, 1} - \sqrt{3}\h{0.5pt}u_{\nu, 2}}{\sqrt{\left(u_{\nu, 1} - \sqrt{3}\h{0.5pt}u_{\nu, 2}\right)^2 + 4u_{\nu, 3}^2}} \right)^{-1/2} \mathbf{e}_{z}.
\end{split}
\end{equation}
Moreover, the director field satisfies the following limit for any $z_0 \in [1, 1 + \epsilon)$:
\begin{align}\label{lim of dir}
\mathbf{d} \left[u_\nu\right] \to \mathbf{e}_{\rho},\h{10pt}\text{as $(\rho, z) \to (0, z_0)$ with $(\rho, z) \in \overline{\mathbb D}^{\h{0.5pt}c} \cap \Big\{ (\rho, z) : \rho^2 + (z - 1)^2 < \epsilon^2 \Big\}$},
\end{align}which identifies the disclinations at the poles as boojum-type surface defects.
\end{enumerate}
\end{armathm}
\begin{armarem}We would like point out that
\begin{enumerate}
\item[$\mathrm{(1).}$] By Item (2) of Theorem \ref{main theorem 1} and the far-field condition \eqref{far field f2}, for each $\nu > \nu_0$, there must be an odd number of singularities located above and below the colloid $B_1$, respectively. The local structure of these singularities is detailed in \cite{yu2020disclinations}. \vspace{0.2pc}
\item[$\mathrm{(2).}$] It follows from \eqref{lim of dir} that two segments on $l_z$ emanate from the poles $\mathcal{N}$ and $\mathcal{S}$, along which the director field exhibits singularities of strength $1$.
\end{enumerate}
\end{armarem}

We utilize partial regularity results to establish Item (1) in Theorem~\ref{main theorem 1}, which hinges on proving the Hölder continuity of the solution at all expected regular points. Away from the \(z\)-axis, the argument is straightforward; by excluding a cylindrical neighborhood of the axis, the singular weight \(\rho ^{-1}\) in the energy density remains uniformly bounded, thereby reducing the problem to a standard two-dimensional framework. Conversely, in neighborhoods of the north and south poles, this weight becomes severely singular—it is unbounded and, in fact, fails to be integrable. To overcome this difficulty, we follow the classical approach of Schoen and Uhlenbeck \cite{schoen1983boundary}: in Sections \ref{cha of va}-\ref{prf of 2.3}, we introduce a boundary-flattening transformation. This coordinate change allows us to establish the desired boundary regularity by leveraging a small-energy estimate coupled with a suitable blow-up argument near the flattened boundary. Equipped with the boundary regularity near \(\mathcal{N}\) and \(\mathcal{S}\), we investigate the local structure of boojums for suitably large \(\nu \) in Section \ref{booj}. Item (2) of Theorem \ref{main theorem 1} is subsequently established in Section \ref{(2) of main 1} via a blow-up argument centered on a family of shrinking half-balls, \(B_{\nu ^{-1}}^{+}\), as \(\nu \to \infty\). Items (3) and (4) are proved in Section \ref{ls of boo}, where the core argument relies on Serrin's maximum principle to deduce the strict positivity of \(u_{\nu ,1}\) in \(\overline{\mathbb{D}}^{c}\). Ultimately, we obtain Item (3) by this strict positivity and the regularity of the solution near the poles.\vspace{0.2pc}

Our next result concerns the coexistence of Saturn ring disclinations within the solutions. To establish this, we impose additional symmetry requirements on the solution space. More precisely, we define  \begin{align*}\mathcal{F}_{\nu ,\mathcal{R}}:=\Big\{w\in \mathcal{F}_{\nu }:w\text{\ is\ }\mathcal{R}\text{-symmetric}\Big\}.\end{align*} In addition to the conditions prescribed for $u^{\circ }$ in \eqref{assump 1}, we assume that  \begin{align}\label{assum fo u3}\big(u^{\circ}\big)_{1}-&\sqrt{3}\h{0.5pt}\big(u^{\circ}\big)_{2}>0 \h{5pt}\text{ at $(1,0)$}, \qquad  \big(u^{\circ}\big)_{3} \h{1pt} \big\vert_{\partial^+ B^+_{1}} \geqslant 0, \qquad \big(u^{\circ}\big)_3 \text{ is not constant on } \partial^+ B^+_1.\end{align} Here, $\p^+ B^+_1$ denotes the curved boundary of the upper half-ball $B_1^+$. Under these assumptions, we state our main results regarding the Saturn ring disclination and its local structural properties. Note that in what follows, $D_r\left(x\right)$ denotes the disk in the $(\rho,z)$-plane with center $x$ and radius $r$.

\begin{armathm}\label{main theorem 3}
Suppose that \(w^{\circ}=L[u^{\circ}]\) is \(\mathcal{R}\)-symmetric, and let \(w_{\nu } = L[u_\nu]\) be a minimizer of \(F_{\nu }\) within the configuration space \(\mathcal{F}_{\nu ,\mathcal{R}}\). In addition, assume that \(u^{\circ }\) satisfies both \eqref{assump 1} and \eqref{assum fo u3}. Then, there exists a \(\nu_{0}>0\) such that for all \(\nu>\nu_{0}\), \(w_{\nu }\) possesses a Saturn ring disclination located at some \((r_{\nu},0)\) with \(r_{\nu}>1\). More precisely, the corresponding \(u_{\nu}=(u_{\nu,1},u_{\nu,2},u_{\nu,3})^{\top}\) satisfies the following properties:
\begin{enumerate}[label=\textup{\arabic*.}, leftmargin=*]
    \item[$\mathrm{(1).}$] There exists a sufficiently small constant $\epsilon = \epsilon_{\nu}>0$ such that $$u_{\nu, 1} - \sqrt{3}\h{0.5pt}u_{\nu,2} > 0  \h{15pt}\text{ on $\big\{(\rho, 0): \rho \in\left(r_{\nu}-\epsilon, r_{\nu}\right)\big\}$} $$ and $$ u_{\nu, 1} - \sqrt{3}\h{0.5pt}u_{\nu,2} < 0  \h{15pt}\text{on $\big\{(\rho, 0): \rho \in\left(r_{\nu}, r_{\nu}+\epsilon\right)\big\}$}.$$ In addition, $u_{\nu}$ is biaxial on $D_{\epsilon}\left(r_{\nu}, 0\right) \setminus\big\{\left(r_{\nu}, 0\right)\big\}$.  More precisely, the three eigenvalues in \eqref{eigenvalues} computed in terms of the components of $u_{\nu}$ satisfy   $$\lambda_3>\lambda_2>\lambda_1 \h{15pt}\text{ on $D_{\epsilon}\left(r_{\nu}, 0\right) \setminus\big\{\left(r_{\nu}, 0\right)\big\}$. }$$

\item[$\mathrm{(2).}$]  The director field $\mathbf{d}$ is continuous on $D_{\epsilon}(r_{\nu}, 0) \setminus \big\{(\rho, 0) : \rho \in [r_{\nu}, r_{\nu}+\epsilon)\big\}$, but discontinuous along the line segment $\big\{(\rho, 0) : \rho \in (r_{\nu}, r_\nu +\epsilon)\big\}$. For any fixed $\epsilon^{\prime} \in (0, \epsilon)$, approaching $(r_{\nu}+\epsilon^{\prime}, 0)$ along $\partial D_{\epsilon^{\prime}}(r_{\nu}, 0)$ yields limiting values of $-\mathbf{e}_{z}$ from the lower half and $\mathbf{e}_{z}$ from the upper half. As we loop clockwise around $\partial D_{\epsilon^{\prime}}(r_{\nu}, 0)$ back to the starting point, $\mathbf{d}$ varies continuously from $-\mathbf{e}_{z}$ to $\mathbf{e}_{z}$, remaining a linear combination of $\mathbf{e}_\rho$ and $\mathbf{e}_z$ with a strictly positive $\mathbf{e}_\rho$-coefficient throughout. This continuous transition from $-\mathbf{e}_{z}$ to $\mathbf{e}_{z}$ characterizes a ring disclination of strength $-\frac{1}{2}$ at $(r_\nu, 0)$.\vspace{0.4pc}

\item[$\mathrm{(3).}$]  Let $\alpha\in [0, 2\pi]$ be an angular variable. The value of $u_{\nu}\left(\left(r_{\nu}, 0\right)+\epsilon^{\prime}\left(\cos \alpha, \sin \alpha\right)\right)$ at $\alpha = 0$ ($2\pi$ resp.) is defined as the limit of $u_{\nu}$ at $(r_{\nu}+\epsilon^{\prime},0)$ along the upper (lower resp.) part of $\partial D_{\epsilon^{\prime}}\left(r_{\nu}, 0\right)$. Then, the tangent map of $\mathbf{d}$, defined as the limit of $\mathbf{d}\left(\left(r_{\nu}, 0\right)+\epsilon^{\prime}\left(\cos \alpha, \sin \alpha\right)\right)$ as $\epsilon^{\prime} \rightarrow 0^{+}$, equals
\[
\begin{cases}
\mathbf{e}_{z}, &  \alpha=0 ; \\[1.5mm]
\frac{\sqrt{2}}{2}\left\{1+\frac{\varkappa \cot \alpha}{\sqrt{4+\varkappa^2 \cot^2 \alpha}}\right\}^{1 / 2} \!\mathbf{e}_{\rho}+\sqrt{\frac{2}{4+\varkappa^2 \cot^2 \alpha}}\left\{1+\frac{\varkappa \cot \alpha}{\sqrt{4+\varkappa^2 \cot^2 \alpha}}\right\}^{-1 / 2} \!\mathbf{e}_{z}, &  \alpha \in(0, \pi) ; \\[1.5mm]
\mathbf{e}_{\rho}, &  \alpha=\pi; \\[1.5mm]
\frac{\sqrt{2}}{2}\left\{1-\frac{\varkappa \cot \alpha}{\sqrt{4+\varkappa^2 \cot^2 \alpha}}\right\}^{1 / 2} \! \mathbf{e}_{\rho}-\sqrt{\frac{2}{4+\varkappa^2 \cot^2 \alpha}}\left\{1-\frac{\varkappa \cot \alpha}{\sqrt{4+\varkappa^2 \cot^2 \alpha}}\right\}^{-1 / 2} \!\mathbf{e}_{z}, &  \alpha \in(\pi, 2\pi) ; \\[1.5mm]
-\mathbf{e}_{z}, &  \alpha=2\pi.
\end{cases}
\]
Here, $\varkappa$ is a non-positive constant given by $\displaystyle \varkappa :=\frac{\partial_\rho u_{\nu,1}-\sqrt{3} \partial_\rho u_{\nu,2}}{\partial_z u_{\nu,3}}\h{1pt}\bigg|_{\left(r_{\nu}, 0\right)}$.
\end{enumerate}
\end{armathm}

The proof of this theorem, provided in Section \ref{saturn}, relies on the first condition in \eqref{assum fo u3} alongside the far-field asymptotic behavior of $w$. This asymptotic behavior is established in Section \ref{far field arg} via a Bochner-type inequality.

\section{Regularity at the Poles}\label{cha of va}

Since the arguments at the north and south poles are identical, we restrict our attention to the north pole without loss of generality. The main result of this section is stated as follows:
\begin{armaprop}\label{main theorem in section 2}
For any $\nu > 0$, there exists a radius $R^*_\nu \in (0,1)$, depending only on $\nu$, such that any minimizer $w_{\nu}$ of $F_{\nu}$ over $\mathcal{F}_{\nu}$ is Hölder continuous on the closure of $B_{R^*_\nu}(\mathcal{N})\cap \Omega$ with some exponent $\alpha\in (0,1)$. Here, $\alpha$ is independent of $\nu$.
\end{armaprop}
Throughout this section, we drop the subscription $\nu$ from the notation whenever no ambiguity can arise. Since $\partial\Omega$ is smooth, we can flatten the part of the boundary near $\mathcal{N}$ using a diffeomorphism. More precisely, the following properties can be satisfied by a suitably chosen diffeomorphism from an upper half-ball. \begin{armalem}\label{properties of coordinates transformation}
For some sufficiently small $\rho_{0}>0$, there is an axially symmetric neighborhood, denoted by $O_{\rho_{0}}(\mathcal{N})$, of $\mathcal{N}$ and a map \begin{equation}\label{diffeomorphism psi}
\psi: \overline{B_{\rho_{0}}^{+}} \coloneq \Big\{(s_{1},s_{2},s_{3}):s_{1}^{2}+s_{2}^{2}+s_{3}^{2}\leqslant\rho_{0}^{2}\h{4pt}\text{and}\h{4pt} s_{3}\geqslant 0\h{1pt}\Big\}\h{1.5pt}\longrightarrow \h{1.5pt} \overline{O_{\rho_{0}}(\mathcal{N})\cap \Omega}
\end{equation} such that $\psi$ is a diffeomorphism to $O_{\rho_{0}}(\mathcal{N})\cap \Omega$ when restricted to $B_{\rho_{0}}^{+}$. In addition, $\psi$ satisfies all properties organized in the following four items.\begin{enumerate}
\item[\textup{(1).}]  The following four facts are held by $\psi$ on $\overline{B_{\rho_{0}}^{+}}$$\mathrm{:}$\vspace{0.4pc}
\begin{enumerate}
     \item[$\mathrm{(1.1).}$] $\psi$ preserves the azimuthal angle. \vspace{0.2pc}
    \item[$\mathrm{(1.2).}$] $\psi(0) = \mathcal{N}$.
    \item[$\mathrm{(1.3).}$] $\psi$ maps the flat boundary $\overline{D_{\rho_{0}}}$ onto $\overline{O_{\rho_{0}}(\mathcal{N})\cap \partial \Omega}$, where $D_{r} \coloneq \Big\{(s_{1},s_{2},0):s_{1}^{2}+s_{2}^{2} < r^{2}\Big\}$.
    \item[$\mathrm{(1.4).}$] The set $\psi\big(\overline{B_\rho^+}\big)$ contains $\overline{B_{\frac{\rho}{2}}(\mathcal{N})\cap \Omega}$ \h{1pt}  for any $\rho\in (0,\rho_{0})$.
\end{enumerate}\vspace{0.6pc}
\item[\textup{(2).}] Let $g$ be the pullback of the Euclidean metric under $\psi$. It is represented by the matrix $(g_{\alpha\beta})$. Its inverse is denoted by $(g^{\alpha\beta})$. Then there are two positive universal constants $C_1$ and $C_2$ such that the following estimates hold for any $\rho \in (0, \rho_{0})$ and $\alpha, \beta \in \{1, 2, 3\}$$\mathrm{:}$
     \begin{align*}
    \big\vert g^{\alpha\beta}-\delta^{\alpha\beta} \big\vert\leqslant C_{1}\rho\quad \h{2pt}\text{and}\h{2pt}\quad
    \left| \sqrt{\det (g_{\alpha\beta})} - 1 \h{1pt} \right| \leqslant  C_{2}\rho  \h{15pt}\text{on $B_\rho^+$.}
    \end{align*}Here, $\delta^{\alpha\beta}$ and $\delta_{\alpha\beta}$ interchangeably denote the Kronecker delta.\vspace{0.6pc}
\item[\textup{(3).}] Let $h=(h_{ij})_{2\times 2}$ be the induced metric matrix corresponding to the parametrization$\mathrm{:}$ $$(s_{1},s_{2}) \in D_{\rho_{0}} \longrightarrow\psi(s_{1},s_{2},0)\in \partial \Omega.$$Here, we still use $D_\rho$ to denote the plane disk $\big\{(s_1, s_2) : s_1^2 + s_2^2 < \rho^2\big\}$.  Then $h$ is smooth on $D_{\rho_{0}}$ and equals the identity matrix at the origin. Moreover, there is another positive universal constant $C_3$ such that the following estimates hold for any $\rho \in (0, \rho_{0})$ and $i, j \in \{1, 2\}$$\mathrm{:}$
\begin{align*}
\big|h_{ij} -\delta_{ij}\big|\leqslant C_{3}\rho\quad\h{2pt}\text{and}\h{2pt}\quad \left| \sqrt{E G -F^{2} } - 1 \h{1pt} \right| \leqslant  C_{3}\rho \h{15pt}\text{on $D_\rho$.}
\end{align*}$E,F,G$ are the coefficients of the first fundamental form associated with the parametrization. \vspace{0.6pc}
\item[\textup{(4).}] Keep taking $\rho_{0}>0$ sufficiently small. Then for any $B_{r}(s)\subset B_{\rho_0}^{+}$, we have \begin{align*}
\mathrm{diam}\,\psi(B_{r}(s))\leqslant \frac{201}{100}r,\quad \mathrm{dist}\big(\psi(s),\partial (\psi(B_{r}(s)))\big)\geqslant \frac{99}{100}r.
\end{align*}
Denote by $\phi$ the inverse of $\psi$ in \eqref{diffeomorphism psi}. For any $B_{r}(x)\subset O_{\rho_0}(\mathcal{N})\cap \Omega$, we also have
\begin{align*}
\mathrm{diam}\,\phi(B_{r}(x))\leqslant \frac{201}{100}r,\quad \mathrm{dist}\big(\phi(x),\partial (\phi(B_{r}(x)))\big)\geqslant \frac{99}{100}r.
\end{align*}
\end{enumerate}

\end{armalem}

Fix $r \in (0, \rho_0)$.  For any $\tilde{v} \in H^1(B_r^+)$, we define \begin{equation}\label{mathfrak F energy}
\mathfrak{F}_{\nu}[\tilde{v},{B_{r}^{+}}] := \int_{B_{r}^{+}} g^{\alpha\beta}  \partial_{s_{\alpha}} \tilde{v} \cdot \partial_{s_{\beta}} \tilde{v} +\sqrt{2}\mu \big(1-3S[\tilde{v}]\big)  \diff V_{g}+\nu\int_{D_{r}}\big\vert \tilde{v}-\tilde{w}^{\circ}\big\vert^{2}\diff S_{g}.
\end{equation}Here, $\tilde{w}^{\circ} := w^\circ \circ \psi$. The volume and surface elements induced by $g$ are respectively given by
    \begin{align*}
\diff V_{g}=\sqrt{\det( g _{\alpha\beta})}\diff s\quad\h{2pt}\text{and}\h{2pt}\quad \diff S_{g}=\sqrt{EG-F^{2}}\diff s_{1}\diff s_{2}.
\end{align*}Let $\p^+ B_r^+$ be the curved boundary of $B_r^+$. A standard filling-hole argument shows that $\tilde{w} := w \circ \psi$ is a minimizer of the energy $\mathfrak{F}_{\nu}[\cdot,{B_{r}^{+}}]$ among all $\tilde{v} \in H_{\mathrm{asy}}^1\big(B_r^+; \mathbb{S}^4\big)$ that satisfy  $\tilde{v} = \tilde{w}$ on $\p^+ B_r^+$. Here, $H_{\mathrm{asy}}^1\big(B_r^+; \mathbb{S}^4\big)$ denotes the collection of maps in $H^1\big(B_r^+; \mathbb{S}^4\big)$ with axial symmetry. The proof of Proposition \ref{main theorem in section 2} is then reduced to proving the H\"{o}lder regularity of $\tilde{w}$ at $0$, which can be established using the following lemma.

\begin{armalem}[Energy decay estimates]\label{Campanato Criteria}\, There is an exponent $\alpha \in (0, 1)$ such that for any  $\nu>0$, there exist a radius $\rho_{*} < \rho_0$ and a constant $\mathscr{C} > 0$, with which the following four energy estimates hold $\mathrm{:}$
\begin{enumerate}
    \item[\textup{(1).}] For any $ r\in \big(0,2^{-3}\rho_{*}\big]$, it satisfies \begin{align*}  \int_{B_{r}^{+}}\vert \nabla \tilde{w} \vert^{2}\leqslant \mathscr{C}r^{1 + 2\alpha}.\end{align*}
    \item[\textup{(2).}] For any $P \in D_{2^{-7}\rho_{*}}$, if it satisfies $0 < 19r< | P |$, then we have
    \begin{align*}  \int_{B_{8r}^{+}(P)}\vert \nabla \tilde{w}\vert^{2}\leqslant \mathscr{C} r^{1 + 2\alpha}.\end{align*}
    \item[\textup{(3).}] For any $s = (0, 0, s_3)$, if it satisfies $0< 19r < s_3 < 2^{-7}\rho_{*}$, then we have \begin{align*}  \int_{B_{8r}(s)}\vert \nabla \tilde{w} \vert^{2}\leqslant \mathscr{C} r^{1 + 2\alpha}.\end{align*}
    \item[\textup{(4).}] For any $s = (s_1, s_2, s_3) \in B^+_{2^{-7}\rho_{*}}$ and $ r\in \big(0,2^{-9}\rho_{*}\big)$, if \h{1pt}$\min\big\{s_3, \sqrt{s_1^2 + s_2^2}\h{2pt}\big\}>3r$, then we have \begin{align*}  \int_{B_{\frac{3r}{2}}(s)}\vert \nabla \tilde{w}\vert^{2}\leqslant \mathscr{C} r^{1 + 2\alpha}.\end{align*}
\end{enumerate}The constants $\rho_*$ and $\mathscr{C}$ only depend on $\nu$ and $\rho_0$.
\end{armalem}

Lemma \ref{Campanato Criteria} will be proved in Sections \ref{sma en} and \ref{prf of 2.3}. Now, we prove Proposition \ref{main theorem in section 2} using this lemma.

\begin{proof}[\textit{Proof of Proposition \ref{main theorem in section 2}}]
Fix a $r' \in \left(0, 2^{-7} \rho_*\right)$, where $\rho_*$ is given in Lemma \ref{Campanato Criteria},  and let $\mathscr{I}$ be the family of balls $\big\{B_{r}(s) : s \in B^{+}_{r'} \h{4pt}\text{and}\h{4pt} 0 < r < 2r'\big\}$. For $s \in B_{r'}^{+}$, we define $P_{s,1}$ to be the shortest-distance projection of $s$ to the flat boundary of $B_{r'}^+$. In addition, we let $P_{s, 2} := s - P_{s, 1}$. By $P_{s, 1}$ and $P_{s, 2}$, the family $\mathscr{I}$ can be represented as a union of four mutually disjoint sub-families of balls, which are listed as follows:
\begin{align*}
\mathscr{I}_{1} &\coloneq \big\{B_{r}(s) \in \mathscr{I} : |P_{s, 1}| \leqslant 4r, \, |P_{s, 2}| \leqslant 4r \big\}, \h{8pt}
\mathscr{I}_{2} \coloneq \big\{B_{r}(s) \in \mathscr{I} : |P_{s, 1}| > 4r, \, | P_{s, 2} | \leqslant 4r \big\}, \\[1.5mm]
\mathscr{I}_{3} &\coloneq \big\{B_{r}(s) \in \mathscr{I} : | P_{s, 1} | \leqslant 4r, \, |P_{s, 2}| > 4r \big\}, \h{10pt}
\mathscr{I}_{4} \coloneq \big\{B_{r}(s) \in \mathscr{I} : |P_{s, 1}| > 4r, \, |P_{s, 2}| > 4r \big\}.
\end{align*}
We further partition $\mathscr{I}_2$ by writing $\mathscr{I}_{2} =  \mathscr{I}_{2,1}\cup \mathscr{I}_{2,2}$, where
\begin{align*}
\mathscr{I}_{2,1} \coloneq \big\{B_{r}(s)\in \mathscr{I}_{2}: \vert P_{s, 1}\vert\leqslant 19r\big\},\h{8pt} \mathscr{I}_{2,2} \coloneq \big\{B_{r}(s)\in \mathscr{I}_{2}: \vert P_{s, 1}\vert> 19r\big\}.\end{align*}
\noindent Similarly, $\mathscr{I}_3$ can be  decomposed into $\mathscr{I}_{3} =  \mathscr{I}_{3,1}\cup \mathscr{I}_{3,2}$ with \begin{align*}
\mathscr{I}_{3,1} \coloneq \big\{B_{r}(s)\in \mathscr{I}_{3}: \vert P_{s, 2}\vert\leqslant 19r\big\},\h{8pt} \mathscr{I}_{3,2}  \coloneq \big\{B_{r}(s)\in \mathscr{I}_{3}: \vert P_{s, 2} \vert> 19r\big\}.
\end{align*}

We estimate the Campanato seminorm of $\tilde{w}$ using Poincar\'{e} inequality. Our aim is to verify \begin{align}\label{cmp est}
    \int_{B^{+}_{r'} \h{.8pt}\cap \h{.8pt} B_{r}(s)}\bigg\vert\tilde{w}-\fint_{B^{+}_{r'}\h{0.8pt}\cap \h{0.8pt} B_{r}(s)}\tilde{w}\bigg\vert^{2}\lesssim \mathscr{C} r^{3 + 2 \alpha} \h{15pt}\text{for any $B_r(s) \in \mathscr{I}$.}
\end{align}The cases discussed in what follows are labeled according to the subscripts of the above sub-families.
\begin{enumerate}
    \item[(1).] Suppose $B_{r}(s)\in \mathscr{I}_{1}$. It holds, for any $s'\in B_{r}(s)$, that $\vert s'\vert\leqslant \left(4\sqrt{2}+1\right)r< 8r$. Hence, $B_{r}(s)\subset B_{8r}$. Using Poincar\'{e} inequality then induces
    \begin{align*}
        \int_{B^{+}_{r'} \h{.8pt}\cap \h{.8pt} B_{r}(s)}\bigg\vert\tilde{w}-\fint_{B^{+}_{r'}\h{0.8pt}\cap \h{0.8pt} B_{r}(s)}\tilde{w}\bigg\vert^{2}\leqslant   \int_{B^{+}_{r'} \h{0.8pt}\cap \h{0.8pt} B_{r}(s)}\bigg\vert\tilde{w}-\fint_{B_{8r}^{+}}\tilde{w}\bigg\vert^{2}
\leqslant  \int_{B_{8r}^{+}}\bigg\vert \tilde{w}-\fint_{B_{8r}^{+}}\tilde{w}\bigg\vert^{2}
\lesssim \h{0.5pt} r^{2}\int_{B_{8r}^{+}}\vert \nabla \tilde{w}\vert^{2}.
\end{align*}Here, $\displaystyle \fint_S \tilde{w}$ denotes the average of $\tilde{w}$ over a set $S$ in $\mathbb R^3$. \eqref{cmp est} follows by (1) in Lemma \ref{Campanato Criteria}. \vspace{0.2pc}
\item[(2.1).] Suppose $B_{r}(s)\in \mathscr{I}_{2,1}$. It holds, for any $s'\in B_{r}(s)$, that $\vert s'\vert\leqslant \left(\sqrt{4^{2}+19^{2}}+1\right)r< 21r$. Note that $4 r < \vert P_{s, 1}\vert  \leqslant \vert s\vert<r'$. We then obtain $B_{r}(s)\subset B_{21r} \subset B_{6 r'}$. \eqref{cmp est} follows using the same arguments above for the sub-family $\mathscr{I}_1$.\vspace{0.4pc}
\item[(2.2).] Suppose $B_{r}(s)\in \mathscr{I}_{2,2}$. It holds, for any $s'\in B_{r}(s)$, that $\vert s' - P_{s, 1} \vert \leqslant \vert s' - s \vert + \vert  P_{s, 2} \vert < 8r$. Note that $19r < \vert P_{s,1}\vert \leqslant \vert s \vert<r'$. We then obtain $B^+_{r}(s)\subset B_{8r}^{+}(P_{s, 1}) \subset B_{2 r'}^+$. The same arguments for $\mathscr{I}_1$ induce that
\begin{align*}
\int_{B^{+}_{r'}\h{0.8pt}\cap \h{0.8pt} B_{r}(s)}\bigg\vert\tilde{w}-\fint_{B^{+}_{r'}\h{0.8pt}\cap\h{0.8pt} B_{r}(s)}\tilde{w}\bigg\vert^{2}
\lesssim  r^{2}\int_{B_{8r}^{+}(P_{s,1})}\vert \nabla \tilde{w}\vert^{2}.
\end{align*}\eqref{cmp est} follows by applying (2) in Lemma \ref{Campanato Criteria} to the right-hand side above.\vspace{0.4pc}
\item[(3.1).] Suppose $B_{r}(s)\in \mathscr{I}_{3,1}$. It holds, for any $s'\in B_{r}(s)$, that $\vert s'\vert\leqslant \left(\sqrt{4^{2}+19^{2}}+1\right)r< 21r$. Note that $4r < \vert P_{s, 2}\vert \leqslant \vert s \vert<r'$. We then obtain $ B_{r}(s)\subset B_{21r} \subset B_{6r'}$. \eqref{cmp est} follows using the same arguments above for the sub-family $\mathscr{I}_1$. \vspace{0.4pc}
\item[(3.2).] Suppose $B_{r}(s)\in \mathscr{I}_{3,2}$. It holds, for any $s'\in B_{r}(s)$, that $\vert s' - P_{s, 2} \vert \leqslant \vert s' - s \vert + \vert  P_{s, 1} \vert < 8r$. Note that $19r < \vert P_{s,2}\vert \leqslant \vert s \vert<r'$. We then obtain $B^+_{r}(s)\subset B_{8r}(P_{s, 2}) \subset B_{2 r'}^+$. Therefore,
\begin{align*}
\int_{B^{+}_{r'}\h{0.8pt}\cap\h{0.8pt} B_{r}(s)}\bigg\vert\tilde{w}-\fint_{B^{+}_{r'}\h{0.8pt}\cap \h{0.8pt} B_{r}(s)}\tilde{w}\bigg\vert^{2}
\lesssim  r^{2}\int_{B_{8r}(P_{s,2})}\vert \nabla \tilde{w}\vert^{2}.
\end{align*}\eqref{cmp est} follows by applying (3) in Lemma \ref{Campanato Criteria} to the right-hand side above.\vspace{0.4pc}
\item[(4).] Suppose $B_{r}(s)\in \mathscr{I}_{4}$. It holds $B_{\frac{3r}{2}}(s) \subset B_{2r'}^+$. Here, we use $4\sqrt{2}r < \vert s \vert<r'$. Therefore,
\begin{align*}
\int_{B^{+}_{r'}\h{0.8pt}\cap\h{0.8pt} B_{r}(s)}\bigg\vert\tilde{w}-\fint_{B^{+}_{r'}\h{0.8pt}\cap\h{0.8pt} B_{r}(s)}\tilde{w}\bigg\vert^{2}
\lesssim  r^{2}\int_{B_{r}(s)}\vert \nabla \tilde{w}\vert^{2}\lesssim \mathscr{C} r^{3 + 2\alpha}.
\end{align*}The second estimate follows from (4) in Lemma \ref{Campanato Criteria}.
\end{enumerate}

The proof is complete by the equivalence between the Morrey-Campanato space and the H\"{o}lder space. See Theorem 5.5 in \cite{giaquinta2013introduction}.
\end{proof}

\section{Smallness of Scaled Energy}\label{sma en}

We define the scaled energy of $\tilde{w}$ on the half-ball $B_r^+$ by
\begin{align*}
\Theta \left(\tilde{w}, 0, r\right)\coloneq \frac{1}{r}\int_{B_{r}^{+}}\vert \nabla \tilde{w}\vert^{2},\quad 0<r<\rho_{0}.
\end{align*}
The main result in this section is the smallness of the scaled energy, as shown below. \begin{armaprop}\label{case 1, Uniformly smallness of the energy density}\,
For any $\varepsilon>0$, there is $R_{\varepsilon} \in (0, \rho_0)$ such that $$\Theta\left(\tilde{w},0,r\right)<\varepsilon\qquad \textup{for any } r\in (0,R_{\varepsilon}].$$ The radius $R_{\varepsilon}$ depends only on $\rho_{0}$, $\nu$ and $\varepsilon$. It tends to $0$ if $\varepsilon$ is taken to $0$.
\end{armaprop}

First, we introduce a new energy $\mathbf{F}_{\nu}[\cdot;B_{r}^{+}]$ to approximate the original energy $\mathfrak{F}_{\nu} [\cdot;B_{r}^{+}]$ (see \eqref{mathfrak F energy}).
\begin{armadef}For any $r \in (0, \rho_0)$ and $\tilde{v} \in H^1(B_r^+; \mathbb{S}^4)$, we define
\begin{align*}
    \mathbf{F}_{\nu}[\tilde{v};B_{r}^{+}]\coloneq \int_{B_{r}^{+}}\vert \nabla \tilde{v}\vert^{2} +\nu\int_{D_{r}} \vert \tilde{v}-\tilde{w}^{\circ}\vert^{2}.
\end{align*}The metrics in the two integrals above are all Euclidean.
\end{armadef}
The following lemma estimates the difference between the approximate and the original energies.
\begin{armalem}\label{difference between the approximation energy and the original energy}
There exists a universal constant $C_{6} > 0$ such that
\begin{align*}
\big|\h{1pt} \mathfrak{F}_{\nu}[\tilde{v};B_{r}^{+}]-\mathbf{F}_{\nu}[\tilde{v};B_{r}^{+}]\h{1pt}\big|\leqslant C_{6} \h{1pt} r \h{1pt}\mathbf{F}_{\nu}[\tilde{v};B_{r}^{+}]+C_{6} \h{1pt} r^{3},
\end{align*}for any $0 < r\leqslant\rho_{0}$ and $\tilde{v} \in H^1\big(B_r^+; \mathbb{S}^4\big)$.
\end{armalem}
\begin{proof} We decompose $\mathfrak{F}_{\nu}[\tilde{v};B_{r}^{+}]-\mathbf{F}_{\nu}[\tilde{v};B_{r}^{+}]$ into four parts:
\begin{align*}
\mathfrak{F}_{\nu}[\tilde{v};B_{r}^{+}]-\mathbf{F}_{\nu}[\tilde{v};B_{r}^{+}]=I_{1,1}+I_{1,2}+I_{2}+I_{3},
\end{align*}
where
\begin{align*}
I_{1,1} &\coloneq \int_{B_{r}^{+}}  \left(g^{\alpha\beta}-\delta^{\alpha\beta}\right) \partial_{s_{\alpha}}\tilde{v}\cdot \partial_{s_{\beta}}\tilde{v} \h{1pt} \sqrt{\mathrm{det}(g_{\alpha\beta})}, \h{20pt}
I_{1,2} \coloneq \int_{B_{r}^{+}}  \vert \nabla \tilde{v} \vert^2 \left(\sqrt{\mathrm{det}(g_{\alpha\beta})}-1 \right), \\[1mm]
I_{2} &\coloneq\sqrt{2}\mu \int_{B_{r}^{+}} \big(1-3S[\tilde{v}] \h{1pt}\big) \sqrt{\mathrm{det}(g_{\alpha\beta})},  \hspace{58pt}
I_{3} \coloneq \nu \int_{D_{r}}  |\tilde{v}-\tilde{w}^{\circ}|^2 \left( \sqrt{EG-F^{2}}-1 \right).
\end{align*}
Apply (2)-(3) in Lemma \ref{properties of coordinates transformation} to estimate these 4 integrals. The proof is complete.
\end{proof}

As a consequence, we obtain the following ``almost minimizing'' property for $\tilde{w}$.
\begin{armalem}\label{energy comparison property}

There exists a universal constant $C_{7} > 0$ such that
\begin{align*}
\mathbf{F}_{\nu}[\tilde{w};B_{r}^{+}]\leqslant \left(1+C_{7}\h{.5pt}r\right)\mathbf{F}_{\nu}[\tilde{v};B_{r}^{+}]+C_{7} \h{1pt}r^{3},
\end{align*}for any $0 < r \leqslant \min\big\{ \h{0.5pt}\rho_0, (2 C_6)^{-1} \big\}$ and $\tilde{v} \in H_{\mathrm{asy}}^1\big(B_r^+; \mathbb{S}^4\big)$ with $\tilde{v} = \tilde{w}$ on $\p^+ B_r^+$.
\end{armalem} The proof is a direct application of Lemma \ref{difference between the approximation energy and the original energy} and the minimality of $\tilde{w}$. We omit it here. Instead, we use this lemma to prove a monotonicity result satisfied by $\tilde{w}$.\begin{armalem}\label{almost monotonicity for type 1}

There exists a universal constant $C_{12}>0$ such that the function
\begin{align*}
g\left(\rho\right)\coloneq \left(1+C_{12}\h{0.5pt}\rho\right)\Theta\left(\tilde{w},0,\rho\right)+C_{12}\left(1+\nu\right)\rho\end{align*}
is non-decreasing with respect to $\rho\in \left(0,\min\big\{\h{0.5pt}\rho_{0}, (2C_{6})^{-1}\big\}\h{0.5pt}\right)$. More precisely, it satisfies \begin{equation}\label{energy comparison 3}
\begin{aligned}
g\left(\rho_1\right) +\int_{\rho_{1}}^{\rho_{2}}\frac{1}{r} \int_{\p^+ B_r^+} \vert \partial_{r}\tilde{w}\vert^{2}   \leqslant  g\left(\rho_2\right), \h{15pt}\text{for any $0 < \rho_1 < \rho_2 < \min\big\{\h{0.5pt}\rho_{0}, (2C_{6})^{-1}\big\}$.}
\end{aligned}
\end{equation}
\end{armalem}
\begin{proof} Denote by $\tau$ the tangential variables along the unit sphere in $\mathbb{R}^3$ with center $0$. Fixing $r \in (0, \rho_0)$ satisfying $\displaystyle \int_{\partial^{+}B^+_{r}}\vert \nabla_{\tau} \tilde{w} \vert^{2} <\infty,$ we introduce a 0-homogeneous comparison map, denoted by $v_{r}$, on $B_{r}^{+}$ as follows:
\begin{align*}
v_r(s)\coloneq
\tilde{w} \big(r \hat{s}\big),\qquad s\in B_{r}^{+} \setminus \{0\}.
\end{align*}
Here, $\hat{s}$ is the normalized vector of $s$. Direct calculations yield
\begin{equation}\label{energy comparison 2}
\begin{aligned}
\int_{B_{r}^{+}}\vert \nabla v_{r}\vert^{2} = r\int_{\partial^{+}B^+_{r}}\vert \nabla_{\tau}\tilde{w}\vert^{2}   = r  \frac{\diff}{\diff r}\int_{B_r^+} \vert\nabla \tilde{w}\vert^{2} - r\int_{\partial^{+}B^+_{r}} \vert \partial_{r} \tilde{w} \vert^{2}.
\end{aligned}
\end{equation}
On the other hand, if we further assume $r < (2C_{6})^{-1}$, then by Lemma \ref{energy comparison property}, it turns out that \begin{align*}
\int_{B_r^+} \vert\nabla \tilde{w}\vert^{2} \leqslant \left(1+C_{12}\h{0.5pt}r\right) \int_{B_r^+}\vert\nabla v_{r}\vert^{2} +C_{12}\left(1+\nu\right)r^{2},
\end{align*}
for some universal constant $C_{12}>0$. Plugging \eqref{energy comparison 2} into the right-hand side above infers
\begin{align*}
\frac{\diff}{\diff r}\left(\frac{1+C_{12}\h{0.5pt}r}{r} \int_{B_r^+} \vert\nabla \tilde{w}\vert^{2} \right)+C_{12}\left(1+\nu\right)\geqslant \frac{1}{r} \int_{\p^+ B_r^+} \vert \partial_{r} \tilde{w}\vert^{2}, \qquad \text{for a.e. } r\in (0,\rho_{0}].
\end{align*}\eqref{energy comparison 3} follows by integrating the variable $r$ in the above inequality from $\rho_1$ to $\rho_2$.\end{proof}

Now, we prove the smallness of $\Theta\left(\tilde{w}, 0 ,r\right)$ for sufficiently small $r$.
\begin{proof}[\textit{Proof of Proposition \ref{case 1, Uniformly smallness of the energy density}}] We divide the proof into 3 steps. \\[2mm]
\noindent \textbf{Step 1.} Suppose to the contrary that Proposition \ref{case 1, Uniformly smallness of the energy density} fails. There exists $\varepsilon_0 > 0$ and a positive sequence $\{r_k\}$ converging to $0$ as $k \to \infty$ such that \begin{align*}\Theta\left(\tilde{w},0,r_k\right) \geqslant \varepsilon_0\qquad \textup{for any } k \in \mathbb{N}.\end{align*}Rescale $\tilde{w}$ and $\tilde{w}^{\circ}$ respectively by letting $$W_{k}(\zeta)\coloneq \tilde{w}(r_{k}\zeta)\quad \text{and}\quad  W^{\circ}_{k}(\zeta)\coloneq \tilde{w}^{\circ}(r_{k}\zeta).$$  It then follows from the contradictory assumption that \begin{align}\label{cont assump}\int_{B_{1}^{+}}\vert \nabla W_{k}\vert^{2}=\frac{1}{r_{k}}\int_{B_{r_{k}}^{+}}\vert \nabla \tilde{w}\vert^{2}=\Theta\left(\tilde{w},0,r_{k}\right) \geqslant \varepsilon_0 \h{15pt}\text{for any $k \in \mathbb{N}$}.
\end{align}By Lemma \ref{almost monotonicity for type 1}, the Dirichlet energy of $W_k$ over $B_1^+$ is uniformly bounded. Up to a subsequence, we can assume $\{W_k\}$ converges to some $W_{\infty}$ as $k \to \infty$, weakly in $H^{1}\left(B_1^+\right)$ and strongly in $L^{2}\left(\p^+B_1^+\right)$. In addition, similar arguments in the proof of \cite[Theorem 5.5]{hardt1992axially} induce that $\{W_k\}$ indeed converges to $W_{\infty}$ strongly in $H^1\left(B_1^+\right)$ as $k \to \infty$.\vspace{0.4pc}

\noindent\textbf{Step 2.} We show the $0$-homogeneity of $W_\infty$ in this step. Since the function $g$ is non-decreasing on the interval $\left(0,\min\big\{\h{0.5pt}\rho_{0}, (2C_{6})^{-1}\big\}\h{0.5pt}\right)$, then by \eqref{cont assump} and Lemma \ref{almost monotonicity for type 1}, there exists a constant $L \geqslant \varepsilon_0$ such that
\begin{align*}
\lim_{k \to \infty}\Big[\left(1+C_{12}\h{0.5pt}r_k\right)\Theta\left(\tilde{w},0, r_k\right)+C_{12}\left(1+\nu\right)r_k\Big]=\lim_{k \to \infty}\Theta\left(\tilde{w},0,r_k\right) = L \geqslant \varepsilon_0.
\end{align*}
Applying the inequality \eqref{energy comparison 3}, we have
\begin{align*}
\int_{0}^{1} \frac{1}{r} \int_{\partial^{+} B^+_{r}} \vert \partial_{r} W_k\vert^{2} = \int_{0}^{r_{k}} \frac{1}{r} \int_{\partial^{+} B^+_{r}} \vert \partial_{r}\tilde{w}\vert^{2} \leqslant g\left(r_{k}\right)-L \h{20pt}\text{for any $k \in \mathbb{N}$.}
\end{align*}In our application of \eqref{energy comparison 3}, we have taken $\rho_1 \to 0$ and let $\rho_2 = r_k$. Using the strong $H^1$-convergence of $\{W_k\}$ obtained in Step 1, we can take $k \to \infty$ in the above estimate and infers  $\partial_{r} W_{\infty}=0$, a.e. in $B^{+}_{1}$. Thus, $W_{\infty}$ is $0$-homogeneous over $B_1^+$. \vspace{0.4pc}

\noindent\textbf{Step 3.}  Now, we prove that $W_{\infty}$ is constant on $B_1^+$, by which a contradiction to \eqref{cont assump} follows. Let $\eta$ be the polar angle, and $\omega$ be the azimuthal angle. Under the axially symmetric ansatz, there is $v=v(\eta)=(v_{1},v_{2},v_{3})^{\top}\in \mathbb{S}^{2}$ such that $W_{\infty}= L[v]$. Here, $\eta\in [0,\frac{\pi}{2}]$ since the domain is $B_1^+$. Note that $W_{\infty}$ is a $0$-homogeneous harmonic map on $B_1^+$. The vector field $v$ must solve the ODE system shown below:
\begin{align}\label{Tangent map ODE}
-(\sin \eta\h{0.5pt} v')' +\frac{1}{\sin \eta}\begin{pmatrix} 4v_{1}\\0\\v_{3}
\end{pmatrix}=\bigg\{\vert v'\vert^{2}\sin \eta+\frac{1}{\sin \eta}\left(4v_{1}^{2}+v_{3}^{2}\right)\bigg\}\h{0.5pt}v,\qquad \eta\in \left(0,\frac{\pi}{2}\right).
\end{align}The notation $'$ refers to the derivative with respect to the $\eta$-variable.

Taking the inner product  with $v'$ on both sides of \eqref{Tangent map ODE} yields
\begin{align*}
     \vert v' \vert^{2}\sin^{2}\eta-\left(4v_{1}^{2}+v_{3}^{2}\right)=C_{e},\h{15pt}\text{for any $\eta\in\big[0,\frac{\pi}{2}\big]$.}
\end{align*}
Here, $C_{e}$ is a constant. Since $W_\infty$ is axially symmetric and has finite Dirichlet energy on $B_1^+$, it turns out $v_1(0) = v_3(0) = 0$. Let $\eta = 0$. The left-hand side above is equal to $0$. Therefore, $C_e = 0$. In addition, using the Robin boundary condition in \eqref{EL of w}, we obtain $\p_z W_\infty \equiv \boldsymbol{0}$ on the flat boundary of $B_1^+$, which equivalently shows that $v'(\frac{\pi}{2}) = \boldsymbol{0}$. Taking  $\eta=\frac{\pi}{2}$ on the left-hand side of the above equality, we then obtain
$\left| v'\right|^{2}(\frac{\pi}{2}) - \left(4v_{1}^{2}\left(\frac{\pi}{2}\right)+v_{3}^{2}\left(\frac{\pi}{2}\right)\right) = 0,$ which forces $v_{1}(\frac{\pi}{2})=v_{3}(\frac{\pi}{2})=0$. Therefore, the only solution of \eqref{Tangent map ODE} is the constant map $v \equiv (0, \pm 1, 0)^{\top}$. The proof is complete.\end{proof}

\section{Proof of Lemma \ref{Campanato Criteria}}\label{prf of 2.3}

In this section, we prove Lemma \ref{Campanato Criteria}. Note that Items (1)-(2) in Lemma \ref{Campanato Criteria} are boundary cases, while Items (3)-(4) are interior cases. Since the proofs for the interior cases are similar to and indeed simpler than the boundary cases, we only show the proofs for Items (1)-(2) in detail.

\subsection{Proof of (1) in Lemma \ref{Campanato Criteria}}

To prove (1) in Lemma \ref{Campanato Criteria}, it suffices to show
\begin{armalem}\label{case 1 lemma}
There exist two constants $0 < \varepsilon_{1}, \lambda_{1} < \frac{1}{2}$, and one universal constant $0 < \theta_{1} < \frac{1}{16}$, such that if $0<r<r_1$, then at least one of the following estimates holds$\mathrm{:}$
\begin{equation*}
\begin{aligned}
(1).\,\,\Theta\left(\tilde{w},0,\theta_1\lambda_1 r\right) \leqslant \frac{1}{2}\,\Theta\left(\tilde{w},0,\lambda_1 r\right),\qquad (2).\,\,
\Theta\left(\tilde{w},0,\theta_1\lambda_1 r\right) \leqslant r^{\frac{1}{2}}.
\end{aligned}
\end{equation*}Here, $r_1$ is $R_{\varepsilon}$ in Proposition \ref{case 1, Uniformly smallness of the energy density} with $\varepsilon = \varepsilon_1$ there.
\end{armalem} \noindent In fact, it can be shown by Lemma \ref{case 1 lemma} that \begin{align*}
\Theta\left(\tilde{w},0,\theta_1\rho\right) + \frac{4}{\sqrt{\lambda_1}} \left(\theta_1 \rho\right)^{\frac{1}{2}} < \frac{1}{2} \left( \Theta\left(\tilde{w},0,\rho\right) + \frac{4}{\sqrt{\lambda_1}} \h{1pt} \rho^{\frac{1}{2}}\right) \h{20pt}\text{for any $\rho \in (0, \lambda_1 r_1)$}.
\end{align*}Recall $g$ in Lemma \ref{almost monotonicity for type 1} and the smallness result in Proposition \ref{case 1, Uniformly smallness of the energy density}. It turns out \begin{align*}
    g\left(\theta_1 \rho\right) + \frac{4}{\sqrt{\lambda_1}} \left(\theta_1 \rho\right)^{\frac{1}{2}} < \frac{1}{2} \left( g\left(\rho\right) + \frac{4}{\sqrt{\lambda_1}} \h{1pt} \rho^{\frac{1}{2}}\right) \h{20pt}\text{for any $\rho \in (0, \lambda_1 r_1)$}.
\end{align*} Since $g$ does not decrease on $\left(0,\min\big\{\h{0.5pt}\rho_{0}, (2C_{6})^{-1}\big\}\h{0.5pt}\right)$, the standard iteration yields \begin{align}\label{a decay es}
    g\left(\rho\right) \leqslant \mathscr{C} \rho^{\alpha_1}\h{20pt}\text{for any $\rho \in \big(0, \min\big\{\h{0.5pt}\rho_{0}, \lambda_1 r_1, (2C_{6})^{-1}\big\} \big)$}.
\end{align}Here, $\mathscr{C} > 0$ is a large constant. $\alpha_1 = - \frac{\ln 2}{\ln \theta_1}$ is a universal constant. Item (1) in Lemma \ref{Campanato Criteria} then follows. \vspace{0.2pc}

The remainder of the section is devoted to the proof of Lemma \ref{case 1 lemma}. We argue by contradiction. Suppose that Lemma \ref{case 1 lemma} fails. Then there exist $\{\lambda_{k}\},\{r_{k}\}$, satisfying $\lambda_{k}\to 0$ and
$\Theta\left(\tilde{w},0,r_{k}\right)\to 0$ as $k\to \infty$, such that both of the following two statements hold for any $k \in \mathbb{N}$:
\begin{align}\label{both and 1}
(1).\,\,\Theta\left(\tilde{w},0,\theta_1 \lambda_{k}r_{k}\right)> \frac{1}{2} \, \Theta\left(\tilde{w},0,\lambda_{k} r_{k}\right),\qquad (2).\,\,\Theta\left(\tilde{w},0, \theta_1 \lambda_{k}r_{k}\right)> r_{k}^{\frac{1}{2}}.
\end{align}
In addition, we define
\begin{align}\label{notat 1}
e_{k}^{2}\coloneq \Theta\left(\tilde{w},0, \lambda_{k} r_{k}\right),\qquad y_{k}\coloneq \fint_{B^{+}_{\lambda_{k}r_{k}}} \tilde{w},\qquad \mathcal{W}_{k}(\zeta)\coloneq \tilde{w}\left(\lambda_{k}r_{k}\zeta\right), \qquad \mathcal{W}_{k}^{rs}\coloneq \frac{\mathcal{W}_{k}-y_{k}}{e_{k}}.
\end{align}
Then, by the Neumann-Poincaré inequality, $\left\{\mathcal{W}^{rs}_{k}\right\}$ is uniformly bounded in $H^{1}(B_{1}^{+})$ with the upper bound depending only on a dimensional constant. Up to a subsequence, we can assume
\[
\mathcal{W}^{rs}_{k} \longrightarrow \mathcal{W}^{rs}_{\infty} \quad
\text{weakly in } H^{1}(B_{1}^{+}) \,\,\text{and  strongly in}\h{3pt} L^{2}(B_{1}^{+}),
\qquad \text{as } k \to \infty.
\]
$\mathcal{W}^{rs}_{\infty}$ is the  limiting map. By the lower semi-continuity, it turns out
\begin{align}\label{lower semi-continuity}
\int_{B_1^+} \left\vert \mathcal{W}_\infty^{rs} \right\vert^2 +  \left\vert \nabla \mathcal{W}_\infty^{rs} \right\vert^2
\;\leqslant\;
\liminf_{k\to\infty}  \int_{B_1^+} \left\vert \mathcal{W}_k^{rs} \right\vert^2 + \left\vert \nabla \mathcal{W}_k^{rs} \right\vert^2.
\end{align}By Fatou's lemma and the $H^1$-regularity of $\mathcal{W}_\infty^{rs}$, there is a $\gamma_1 \in (7/8, 1)$ and a subsequence, which we still denote by $\left\{\mathcal{W}_k^{rs}\right\}$, such that \begin{align}\label{surfa bound}
\sup_{k \h{.8pt}\in \h{.8pt} \mathbb{N}\h{.8pt}\cup \h{.8pt} \{\infty\}} \int_{\p^+ B_{\gamma_1}^+} \left\vert \mathcal{W}_k^{rs} \right\vert^2 + \left\vert \nabla \mathcal{W}_k^{rs} \right\vert^2 = b^2_0 < \infty.
\end{align}Using the compactness of the trace operator, we can also assume $\mathcal{W}_k^{rs} \to \mathcal{W}_\infty^{rs}$ strongly in $L^2(\p^+ B_{\gamma_1}^+)$. Note that $\vert \tilde{w}\vert=1$ and  $e_{k} \to 0$ as $k\to \infty$, it follows $\vert y_{k}\vert \to 1$. \vspace{0.2pc}

Using the above notation, we apply a change of variable and rewrite (1) in \eqref{both and 1} by
\begin{align}\label{inequality to derive the contradiction 1}
\frac{1}{\theta_1}\int_{B_{\theta_1 }^{+}} \left\vert \nabla \mathcal{W}^{rs}_{k} \right\vert^{2}>\frac{1}{2}\qquad \text{for any $k\in \mathbb{N}.$}
\end{align}We will prove that \eqref{inequality to derive the contradiction 1} cannot be satisfied if (2) in \eqref{both and 1} holds. Preliminary lemmas are required before we show the proof. First, we establish the following relationship between $r_{k}$ and $e_{k}$.
\begin{armalem} The sequences $\{r_k\}$ and $\{e_k\}$ defined in \eqref{both and 1}-\eqref{notat 1} satisfy the following limit\h{0.3pt}$\mathrm{:}$
    \label{compare r_{k} and e_{k}}
$$\lim_{k \to \infty}  r_{k} \h{0.3pt}e_{k}^{-2} = 0.$$
\end{armalem}
\vspace{-\baselineskip}
\begin{proof} By the definition of $\Theta$, it holds
\begin{align*}
\left(\theta_1 \lambda_{k}r_{k} \right)  \Theta\left(\tilde{w},0,\theta_1 \lambda_{k}r_{k}\right)\leqslant \left(\lambda_{k}r_{k}\right)  \Theta\left(\tilde{w},0,\lambda_{k}r_{k}\right)=\lambda_{k}r_{k} e_{k}^{2}.
\end{align*}
On the other hand, (2) in $\eqref{both and 1}$ implies
\begin{align*}
\left(\theta_1 \lambda_{k}r_{k}\right)  \Theta\left(\tilde{w},0,\theta_1 \lambda_{k}r_{k}\right)> \left(\theta_1 \lambda_{k}r_{k}\right) r_{k}^{\frac{1}{2}}=\theta_1 \lambda_{k}r_{k}^{\frac{3}{2}}.
\end{align*}
Combining these two estimates yields $\theta_1 \lambda_{k}r_{k}^{\frac{3}{2}}\leqslant \lambda_{k}r_{k} e_{k}^{2}$,
 which infers the limit in this lemma. Note that we also use $e_k \to 0$ as $k \to \infty$.\end{proof}

\begin{armadef} We define a configuration space for the limit of $\big\{\mathcal{W}^{rs}_k\big\}$ as follows:
\begin{align*}
\mathcal{F}_{1}\coloneq \left\{\h{0.5pt}\mathcal{W}\in H^{1}\left(B_{\gamma_1}^{+};\mathbb{R}^{5}\right):\mathcal{W}\,\,\textup{is axially symmetric},\,\, \mathcal{W}=\mathcal{W}_{\infty}^{rs}\,\,\text{on}\,\,\partial^{+}B^+_{\gamma_1}\h{0.5pt}\right\}.
\end{align*}
\end{armadef}Then, we claim that

\begin{armalem}
As $k \to \infty$, $\left\{\mathcal{W}_{k}^{rs}\right\}$ converges to $\mathcal{W}^{rs}_{\infty}$ strongly in $H^{1}\big(B_{\gamma_1}^{+}\big)$. In addition, $\mathcal{W}^{rs}_{\infty}$ minimizes Dirichlet energy on $B_{\gamma_1}^{+}$ in the configuration space $\mathcal{F}_{1}$.
\end{armalem}
\begin{proof} We divide the proof into 3 steps.

\noindent\textbf{Step 1. Comparison map.} Fix $k \in \mathbb{N}$, $R > 0$, $\mathcal{W} \in \mathcal{F}_{1}$ and define $$ \displaystyle
M_{k,R, \mathcal{W}}\coloneq
 y_{k}+Re_{k}\frac{\mathcal{W}}{\vert \mathcal{W} \vert \vee R},\h{20pt}\text{where $a \vee b = \max \big\{ a, b\big\}$ for any $a, b \in \mathbb{R}$.}$$
With $R$ and $\mathcal{W}$ fixed, it follows that \begin{align}\label{unif conv of M}
\left\| M_{k,R, \mathcal{W}} - y_k \right\|_{L^\infty(B_{\gamma_1}^+)} \leqslant R e_k \to 0, \h{20pt}\text{as $k \to \infty$.}
\end{align} Since $|y_k| \to 1$ as $k \to \infty$, the normalized vector field of $M_{k,R, \mathcal{W}}$, that is, $\widehat{M}_{k,R, \mathcal{W}} $, is well-defined for large $k$. The largeness of $k$ depends only on $R$ and is independent of $\mathcal{W}$. Using the above notation and a given $s\in (0,1)$, we define
\begin{equation*} \mathcal{W}_{k,s,R}(\zeta) :=
\begin{cases}
    M_{k,R, \mathcal{W}} \left( \dfrac{\zeta}{1-s} \right), &\h{15pt} \zeta \in B_{\left(1-s\right)\h{0.5pt}\gamma_1}^{+}; \\[4mm]
    \dfrac{\gamma_1 - |\zeta|}{s \gamma_1} \h{1pt} M_{k,R, \mathcal{W}^{rs}_{\infty}} \big( \h{0.8pt} \gamma_1 \widehat{\zeta}\h{1pt} \big) + \dfrac{ |\zeta|- \left(1-s\right) \gamma_1}{s \gamma_1} \h{1pt}\mathcal{W}_{k} \big(\h{0.8pt} \gamma_1 \widehat{\zeta} \h{1pt} \big), &\h{15pt} \zeta \in \overline{B^{+}_{\gamma_1}} \setminus B_{\left(1-s\right)\h{0.5pt}\gamma_1}^{+}.
\end{cases}
\end{equation*}
The normalized vector field of $\mathcal{W}_{k,s,R}$, that is, $\widehat{\mathcal{W}}_{k,s,R}$, is our comparison map. \vspace{0.4pc}

\noindent\textbf{Step 2. Validity of the comparison map.} Using \eqref{surfa bound}, we obtain, by the arithmetic-geometric mean inequality, that \begin{align*} \int_0^{\frac{\pi}{2}} \left\vert\h{1pt} \p_\phi \left( \big(\mathcal{V}_k \big)_1^2 + \big(\mathcal{V}_k\big)_3^2 \right) \h{1.5pt}\right\vert \h{1.5pt}\lesssim\h{1.5pt} \int_0^{\frac{\pi}{2}} \sin \phi \h{1pt}\big\vert \p_\phi \mathcal{V}_k  \big\vert^2 + \frac{\big(\mathcal{V}_k \big)_1^2 + \big(\mathcal{V}_k \big)_3^2}{\sin \phi} \h{1.5pt}\lesssim\h{1.5pt} b_0^2 e_k^2 \h{15pt}\text{at $r = \gamma_1$.}
\end{align*} Here, $\mathcal{W}_k  = L[\mathcal{V}_k]$. $\big(\mathcal{V}_k\big)_j$ denotes the $j$-th component of $\mathcal{V}_k$. Therefore, \begin{align*}
\big(\mathcal{V}_k\big)_1^2 + \big(\mathcal{V}_k\big)_3^2 \leq c^2_* \h{0.3pt}b^2_0 \h{0.3pt} e_k^2 \h{20pt}\text{on $\p^+ B_{\gamma_1}^+$, for all $k \in \mathbb{N}$.}
\end{align*}$c_* > 0$ is a large universal constant. Taking $k$ large, we find $\big| \h{0.5pt} \big(\mathcal{V}_k\big)_2 \h{0.2pt}\big| \geq \frac{1}{2}$ on $\p^+ B_{\gamma_1}^+$. Since \begin{align*}
    \p_\phi \big(\mathcal{V}_k\big)_2 = - \frac{\big(\mathcal{V}_k\big)_1}{\big(\mathcal{V}_k\big)_2} \h{0.5pt}\p_\phi \big(\mathcal{V}_k\big)_1 - \frac{\big(\mathcal{V}_k\big)_3}{\big(\mathcal{V}_k\big)_2} \h{0.5pt}\p_\phi \big(\mathcal{V}_k\big)_3,
\end{align*}it follows that \begin{align*}
    \int_0^{\frac{\pi}{2}} \left\vert\h{1pt} \p_\phi \big(\mathcal{V}_k \big)_2 \h{1.5pt}\right\vert \h{1.5pt}\lesssim\h{1.5pt} b_0^2 e_k^2 \h{15pt}\text{at $r = \gamma_1$.}
\end{align*}Note that $\big(\mathcal{V}_k\big)_2$ must converge to $1$ or $-1$ at some point on $\p^+B_{\gamma_1}^+$. If we suppose that $\big(\mathcal{V}_k\big)_2$ converges to $1$ at some location on $\p^+ B_{\gamma_1}^+$, then the last estimate yields the uniform convergence of $\big(\mathcal{V}_k\big)_2$ to $1$ on $\p^+ B_{\gamma_1}^+$. In this case, we have $\mathcal{V}_k$ uniformly converges to $(0,1,0)^\top$ on $\p^+ B_{\gamma_1}^+$. Using the boundedness of the $L^2$-norm of $\mathcal{W}_k^{rs}$ in \eqref{surfa bound}, we obtain $y_k \to (0,0,1,0,0)^\top$. Therefore, $\widehat{\mathcal{W}}_{k,s,R}$ is well-defined for large $k$ by these convergence results and \eqref{unif conv of M}. The same result holds if $\big(\mathcal{V}_k\big)_2$ converges to $-1$ at some location on $\p^+ B_{\gamma_1}^+$.

We now check the boundary condition of $\widehat{\mathcal{W}}_{k, s, R}$ and its $H^1$-regularity on $B_{\gamma_1}^+$. It suffices to consider $\mathcal{W}_{k, s, R}$. First, for any $\zeta \in \p^+ B_{\gamma_1}^+$, it holds $\mathcal{W}_{k, s, R}(\zeta) = \mathcal{W}_k\big(\gamma_1 \widehat{\zeta}\h{1pt}\big) = \tilde{w}\big(\lambda_k r_k \gamma_1 \widehat{\zeta}\h{1pt}\big)$. Second, for any $\zeta \in \p^+ B_{\left(1 - s\right)\gamma_1}^+$, we have $\mathcal{W}_{k, s,R}(\zeta) = M_{k,R, \mathcal{W}^{rs}_{\infty}} \big( \h{0.8pt} \gamma_1 \widehat{\zeta}\h{1pt} \big)$. Since $\mathcal{W} = \mathcal{W}_\infty^{rs}$ on $\p^+ B_{\gamma_1}^+$, the trace of $M_{k,R, \mathcal{W}} \big(  \zeta/(1-s) \big)$ on $\p^+ B_{\left(1-s\right)\gamma_1}^+$ also equals $M_{k,R, \mathcal{W}^{rs}_{\infty}} \big( \h{0.8pt} \gamma_1 \widehat{\zeta}\h{1pt} \big)$. Thus, $\mathcal{W}_{k, s, R}$ is $H^1$-regular on $B_{\gamma_1}^+$.\vspace{0.2pc}

\noindent\textbf{Step 3. Energy comparison.} Lemma \ref{energy comparison property} implies that
\begin{align}\label{rescaling energy comparison 2}
\liminf_{k\to \infty}\frac{1}{\lambda_{k}r_{k}e_{k}^{2}}\int_{B_{  \lambda_{k}r_{k} \gamma_1}^{+}}\vert \nabla \tilde{w}\vert^{2}
&\leqslant \liminf_{k\to \infty}\frac{1}{\lambda_{k}r_{k}e_{k}^{2}}\mathbf{F}_{\nu}[\tilde{w};B_{ \lambda_{k}r_{k} \gamma_1}^{+}]\\[1mm]
&\leqslant \liminf_{k\to \infty}\frac{1}{\lambda_{k}r_{k}e_{k}^{2}}\left\{\left(1+C_{7}\lambda_{k}r_{k}\right)\mathbf{F}_{\nu}[\tilde{v};B_{ \lambda_{k}r_{k} \gamma_1}^{+}]+C_{7}\left(\lambda_{k}r_{k}\right)^{3}\right\}. \nonumber
\end{align}Here in this step, we specifically take $\tilde{v}(\cdot)\coloneq \widehat{\mathcal{W}}_{k,s,R}\left(\frac{\cdot}{\lambda_{k}r_{k}}\right)$. Applying change of variables and the lower semi-continuity, we get
\begin{align*}
\int_{B_{\gamma_1}^+} \big\vert \nabla \mathcal{W}^{rs}_{\infty} \big\vert^{2} \leqslant \liminf_{k \to \infty}\int_{B_{\gamma_1}^+} \big\vert \nabla \mathcal{W}^{rs}_{k} \big\vert^{2} = \liminf_{k \to \infty}\frac{1}{\lambda_{k}r_{k}e_{k}^{2}}\int_{B_{ \lambda_{k}r_{k} \gamma_1}^{+}}\vert \nabla \tilde{w}\vert^{2}.
\end{align*}
Utilizing Lemma \ref{compare r_{k} and e_{k}} induces
\begin{align*}
&\limsup_{k\to \infty}\frac{1}{\lambda_{k}r_{k}e_{k}^{2}}\int_{D_{ \lambda_{k}r_{k} \gamma_1}}\vert \tilde{v}-\tilde{w}^{\circ}\vert^{2}\lesssim \lim_{k\to \infty}\frac{\lambda_{k}r_{k}}{e^2_{k}}=0\quad \text{and}\quad \lim_{k\to \infty} \lambda_{k}^{2}r_{k}^{2} e_{k}^{-2}=0.
\end{align*}In light of the above estimates and limits, \eqref{rescaling energy comparison 2} can be reduced to \begin{align}\label{reduced comp ine}
\int_{B_{\gamma_1}^+} \big\vert \nabla \mathcal{W}^{rs}_{\infty} \big\vert^{2} \leqslant \liminf_{k\to \infty}\frac{1+C_{7}\lambda_{k}r_{k}}{\lambda_{k}r_{k}e_{k}^{2}} \int_{B_{ \lambda_{k}r_{k} \gamma_1}^{+}} |\h{0.5pt} \nabla \tilde{v} \h{0.8pt}|^2 = \liminf_{k\to \infty}\frac{1+C_{7}\lambda_{k}r_{k}}{e_{k}^{2}} \int_{B_{\gamma_1}^{+}} \big|\h{0.5pt} \nabla \widehat{\mathcal{W}}_{k, s, R} \h{0.8pt}\big|^2.
\end{align}

To calculate the limit on the right-hand side above, we decompose the domain into two parts. First, \begin{align*}
    e_k^{-2}\int_{B_{\left(1 - s \right)\gamma_1}^{+}} \big|\h{0.5pt} \nabla \widehat{\mathcal{W}}_{k, s, R} \h{0.8pt}\big|^2 \leqslant \left(1 - s \right) \int_{B_{\gamma_1}^+} \big| M_{k,R, \mathcal{W}} \big|^{-2} \h{1.5pt}  \left| R  \nabla \frac{\mathcal{W}}{| \mathcal{W} |  \vee R}\right|^2.
\end{align*}Utilizing \eqref{unif conv of M} and the $H^1$-regularity of $\mathcal{W}$, we have \begin{align}\label{limt on big}
\limsup_{s \to 0} \h{1pt}\limsup_{ R \to \infty}\h{1pt} \limsup_{k \to \infty} \h{1pt} e_k^{-2}\int_{B_{\left(1 - s \right)\gamma_1}^{+}} \big|\h{0.5pt} \nabla \widehat{\mathcal{W}}_{k, s, R} \h{0.8pt}\big|^2 \leqslant \int_{B_{\gamma_1}^+}  \big|  \nabla  \mathcal{W} \big|^2.
\end{align} On the shell region $B_{\gamma_1}^+ \setminus B_{\left(1-s\right)\gamma_1}^+$, we estimate the energy of $\widehat{\mathcal{W}}_{k,s,R}$ as follows: \begin{align*}
    &e_k^{-2}\int_{B_{\gamma_1}^+ \setminus B_{\left(1 - s \right)\gamma_1}^{+}} \big|\h{0.5pt} \nabla \widehat{\mathcal{W}}_{k, s, R} \h{0.8pt}\big|^2 \lesssim e_k^{-2}\int_{B_{\gamma_1}^+ \setminus B_{\left(1 - s \right)\gamma_1}^{+}}  \big|\h{0.5pt} \nabla \mathcal{W}_{k, s, R} \h{0.8pt}\big|^2 \\[1.5mm]
    &\h{20pt}\lesssim \int_{\left(1 - s \right)\gamma_1}^{\gamma_1} \int_{\p^+ B_{\gamma_1}^+}  \big| \h{1pt} \nabla \mathcal{W}^{rs}_{\infty}  \h{1pt}\big|^2 +  \big| \nabla \mathcal{W}^{rs}_k\big|^2  + \left(s \gamma_1\right)^{-2}\int_{B_{\gamma_1}^+ \setminus B_{\left(1 - s \right)\gamma_1}^{+}} \left| \mathcal{W}^{rs}_{k}   -   \frac{ R\mathcal{W}_{\infty}^{rs}}{\big| \mathcal{W}_\infty^{rs} \big| \vee R}   \right|^2
\big(\h{0.8pt} \gamma_1 \widehat{\zeta}\h{1pt}\big).  \end{align*}In light of \eqref{surfa bound} and the strong convergence of $\big\{\mathcal{W}_k^{rs}\big\}$ to $\mathcal{W}_\infty^{rs}$ in $L^2(\p^+ B_{\gamma_1}^+)$, it turns out from the last estimate that \begin{align}\label{es on ann}
  \limsup_{s \to 0} \h{1pt}\limsup_{ R \to \infty}\h{1pt} \limsup_{k \to \infty} \h{1pt}  e_k^{-2}\int_{B_{\gamma_1}^+ \setminus B_{\left(1 - s \right)\gamma_1}^{+}} \big|\h{0.5pt} \nabla \widehat{\mathcal{W}}_{k, s, R} \h{0.8pt}\big|^2 = 0.
\end{align}Applying this limit and \eqref{limt on big} to the right-hand side of \eqref{reduced comp ine} infers $\displaystyle
\int_{B_{\gamma_1}^+}\big\vert \nabla \mathcal{W}_{\infty}^{rs} \big\vert^{2} \leqslant \int_{B_{\gamma_1}^+} \big\vert \nabla \mathcal{W} \big\vert^{2}$. The proof is complete.
\end{proof}

Now, we prove (1) in Lemma \ref{Campanato Criteria} by disproving \eqref{inequality to derive the contradiction 1}.
\begin{proof}[\textit{Proof of (1) in Lemma \ref{Campanato Criteria}}] Taking $k \to \infty$ in \eqref{inequality to derive the contradiction 1}, we obtain $\displaystyle \frac{1}{\theta_1}\int_{B_{\theta_1 }^{+}} \left\vert \nabla \mathcal{W}^{rs}_{\infty} \right\vert^{2}  \geqslant \frac{1}{2}.$
In addition, it holds
$\big\Vert \nabla \mathcal{W}_{\infty}^{rs}\big\Vert_{L^{\infty}(B_{3/8}^{+})}\lesssim \big\Vert \nabla \mathcal{W}^{rs}_{\infty}\big\Vert_{L^{2}(B_{\gamma_1}^{+})}\lesssim 1$ by the mean value identity for harmonic functions. In light of these two estimates,  $\theta_1$ must be bounded from below by a positive universal constant, which yields a contradiction if $\theta_1$ is taken suitably small. The smallness of $\theta_1$ is universal.
\end{proof}

\subsection{Proof of (2) in Lemma \ref{Campanato Criteria}}

In this section, we need to use the energy defined on the $(\rho, z)$-plane. That is \begin{align*}
\Xi\left(\tilde{u}, P, r\right)\coloneq \int_{D_{r}^{+}\left(P\right)} \left\vert D \tilde{u} \h{.8pt}\right\vert^{2}+\frac{4\tilde{u}_{1}^{2}+\tilde{u}_{3}^{2}}{\rho^{2}}.
\end{align*}
(2) in Lemma \ref{Campanato Criteria} can then be reduced to the following 2D version: \begin{armalem}\label{2D ver}
There is an exponent $\alpha \in (0, 1)$ such that for any  $\nu>0$, there exist a radius $\rho_{*} < \rho_0$ and a constant $\mathscr{C} > 0$, with which the
energy estimate:
\begin{align*}
   \Xi\left(\tilde{u}, P, 8r\right)
     \leqslant \mathscr{C} r^{\alpha}\end{align*} holds for any $P \in D_{2^{-7}\rho_{*}}$ and $r \in \left(0, \frac{| P |}{19}\right)$.
\end{armalem}
\noindent Indeed, if $P$ satisfies the assumption in Lemma \ref{2D ver}, then $B_{8r}^+\left(P\right) \subset D_{8r}^+\left(P\right) \times \left\{\phi : \left| \phi \right| < \arcsin \frac{8r}{|P|}\right\}$, where $\phi$ is the azimuthal angle. Hence, \begin{align*}
    \int_{B_{8r}^{+}(P)}\vert \nabla \tilde{w}\vert^{2}\leqslant 2 \arcsin \frac{8r}{|P|} \int_{D_{8r}^{+}(P)} \rho \left\vert D \tilde{u} \h{.8pt}\right\vert^{2}+\frac{4\tilde{u}_{1}^{2}+\tilde{u}_{3}^{2}}{\rho}  \leqslant 2 \left( | P | + 8r \right) \arcsin \left(\frac{8r}{|P|} \right) \Xi\left(\tilde{u}, P, 8r\right).
\end{align*}(2) in Lemma \ref{Campanato Criteria} follows from the above estimate and Lemma \ref{2D ver}.\vspace{0.2pc}

To prove Lemma \ref{2D ver}, it suffices to show
\begin{armalem}\label{case 2 lemma 2D}
There exist two constants $0 < \varepsilon_{2}, \lambda_{2} < \frac{1}{2}$, and one universal constant $0 < \theta_{2} < \theta_1 < \frac{1}{16}$, such that if $P \in D_{2^{-7}r_2}$ and $\rho \in \left(0, 2^{-4} | P |\h{0.5pt} \right)$, then at least one of the following estimates holds$\mathrm{:}$
\begin{equation*}
\begin{aligned}
(1).\,\,\Xi\left(\tilde{u}, P,\theta_2\lambda_2 \h{0.5pt}\rho\right) \leqslant \frac{1}{2}\,\Xi\left(\tilde{u}, P,\lambda_2 \h{0.5pt}\rho\right),\qquad (2).\,\,
\Xi\left(\tilde{u},P,\theta_2\lambda_2 \h{0.5pt}\rho\right) \leqslant \rho^{\frac{1}{2}}.
\end{aligned}
\end{equation*}Here, $r_2$ is $R_{\varepsilon}$ in Proposition \ref{case 1, Uniformly smallness of the energy density} with $\varepsilon = \varepsilon_2$ there.
\end{armalem}
\noindent In fact, it can be shown by Lemma \ref{case 2 lemma 2D} that \begin{align*}
\Xi\left(\tilde{u},P,\theta_2 r\right) + \frac{4}{\sqrt{\lambda_2}} \left(\theta_2 r\right)^{\frac{1}{2}} < \frac{1}{2} \left( \Xi\left(\tilde{u},P, r\right) + \frac{4}{\sqrt{\lambda_2}} \h{1pt} r^{\frac{1}{2}}\right) \h{20pt}\text{for any $r \in \big(0, 2^{-4} \lambda_2  | P | \h{0.5pt}\big)$}.
\end{align*}The standard iteration yields \begin{align}\label{tem est}
    \Xi\left(\tilde{u},P, r\right) + \frac{4}{\sqrt{\lambda_2}} \h{1pt} r^{\frac{1}{2}} \leqslant 2 r^{\alpha_2}A^{- \alpha_2} \left( \Xi\left(\tilde{u},P, A\right) + \frac{4}{\sqrt{\lambda_2}} \h{1pt} A^{\frac{1}{2}}\right)\h{20pt}\text{for any $r \in \big(0, A\h{0.5pt}\big)$}.
\end{align}Here, $\alpha_2 = - \frac{\ln 2}{\ln \theta_2}$ is a universal constant. $A$ is a simple notation for $2^{-4} \lambda_2  | P | $. Note that $P \in D_{2^{-7} r_2}$. Hence, $$D_A^+\left(P\right) \times \big\{ \phi : |\h{0.5pt} \phi \h{0.5pt}| < \pi \big\} \subset B^+_{ \left(1 + 2^{-4} \lambda_2 \right) \h{0.5pt}| P |} \subset B^+_{ 2 \h{0.5pt}| P |} \subset B_{r_2}^+,$$ which furthermore infers \begin{align*}
     \pi  | P | \int_{D_{A}^{+}(P)}  \left\vert D \tilde{u} \h{.8pt}\right\vert^{2}+\frac{4\tilde{u}_{1}^{2}+\tilde{u}_{3}^{2}}{\rho^2} \leqslant 2 \pi \int_{D_{A}^{+}(P)} \rho \left\vert D \tilde{u} \h{.8pt}\right\vert^{2}+\frac{4\tilde{u}_{1}^{2}+\tilde{u}_{3}^{2}}{\rho}   \leqslant \int_{B^+_{2 \h{0.5pt} | P |}} \left| \nabla \tilde{w} \right|^2.
\end{align*} It then turns out by \eqref{a decay es} that $\Xi\left(\tilde{u}, P, A\right) \leqslant \Theta \left( \tilde{w}, 0, 2\h{0.5pt}|P| \right) \leqslant \mathscr{C} |P|^{\alpha_1}$. Applying this estimate to \eqref{tem est}, we obtain Lemma \ref{2D ver}. Here, we also use $\theta_2 < \theta_1$. \vspace{0.2pc}

The remainder of the section is devoted to the proof of Lemma \ref{case 2 lemma 2D}. We argue by contradiction. Suppose that Lemma \ref{case 2 lemma 2D} fails. Then there exist $\{\varepsilon_{k}\},\{\lambda_{k}\} \to 0$, $P_k \in D_{2^{-7} R_{\varepsilon_k}}$ and $\rho_k \in \big(0, 2^{-4} | P_k |\h{0.8pt}\big)$ such that both of the following two statements hold for $k \in \mathbb{N}$ suitably large:
\begin{align}\label{both and 2}
(1).\,\,\Xi\left(\tilde{u}, P_k,\theta_2 \lambda_{k} \rho_{k}\right)> \frac{1}{2} \, \Xi\left(\tilde{u},P_k,\lambda_{k} \rho_{k}\right),\qquad (2).\,\,\Xi\left(\tilde{u},P_k, \theta_2 \lambda_{k}\rho_{k}\right)> \rho_{k}^{\frac{1}{2}}.
\end{align}
Here, $R_{\varepsilon_k}$ is given in Proposition \ref{case 1, Uniformly smallness of the energy density} with $\varepsilon = \varepsilon_k$ there. In addition, we define
\begin{align}\label{notat 2}
\tau_{k}^{2}\coloneq \Xi\left(\tilde{u}, P_k, \lambda_{k} \rho_{k}\right),\quad z_{k}\coloneq \fint_{D^{+}_{\lambda_{k} \rho_{k}}\left(P_k\right)} \tilde{u},\quad \mathcal{U}_{k}(\zeta)\coloneq \tilde{u}\left(P_k + \lambda_{k}\rho_{k}\zeta\right), \quad \mathcal{U}_{k}^{rs}\coloneq \frac{\mathcal{U}_{k}-z_{k}}{\tau_{k}}.
\end{align}
By the Neumann-Poincaré inequality, $\big\{\mathcal{U}^{rs}_{k}\big\}$ is uniformly bounded in $H^{1}(D_{1}^{+})$. Using the Rellich-Kondrachov theorem and the Banach-Alaoglu's theorem, we have, up to a subsequence, \begin{align*} \mathcal{U}^{rs}_{k}\longrightarrow \mathcal{U}^{rs}_{\infty} \quad \textup{weakly in } H^{1}(D_{1}^{+})\textup{ and strongly in } L^{2}(D_{1}^{+}).\end{align*}
$\mathcal{U}_{\infty}^{rs} \in H^1(D_1^+)$ is the limiting map. By the lower semi-continuity, it turns out
\begin{align}\label{lower semi-continuity 2D}
\int_{D_1^+} \left\vert \mathcal{U}_\infty^{rs} \right\vert^2 +  \left\vert D \mathcal{U}_\infty^{rs} \right\vert^2
\;\leqslant\;
\liminf_{k\to\infty}  \int_{D_1^+} \left\vert \mathcal{U}_k^{rs} \right\vert^2 + \left\vert D \mathcal{U}_k^{rs} \right\vert^2.
\end{align}By Fatou's lemma and the $H^1$-regularity of $\mathcal{U}_\infty^{rs}$, there is a $\gamma_2 \in (7/8, 1)$ and a subsequence, which we still denote by $\left\{\mathcal{U}_k^{rs}\right\}$, such that \begin{align}\label{cir bound}
\sup_{k \h{.8pt}\in \h{.8pt} \mathbb{N}\h{.8pt}\cup \h{.8pt} \{\infty\}} \int_{\p^+ D_{\gamma_2}^+} \left\vert \mathcal{U}_k^{rs} \right\vert^2 + \left\vert D \mathcal{U}_k^{rs} \right\vert^2 = b^2_0 < \infty.
\end{align}Using Arzel\`{a}-Ascoli theorem, we can also assume $\mathcal{U}_k^{rs} \to \mathcal{U}_\infty^{rs}$ uniformly on $\p^+ D_{\gamma_2}^+$. Note that $\vert \tilde{u}\vert=1$ and  $\tau_{k} \to 0$ as $k\to \infty$, it follows $\vert z_{k}\vert \to 1$. \vspace{0.2pc}

Using the above notation, we apply a change of variables and rewrite (1) in \eqref{both and 2} by
\begin{equation}\label{Contradiction equation in case 2}
\int_{D_{\theta_{2}}^{+}} \left\vert D  \mathcal{U}_{k}^{rs} \right\vert^{2}+\left(\frac{\lambda_{k}\rho_{k}}{\tau_{k}}\right)^{2}\h{1pt}\frac{4 \,\mathcal{U}_{k,1}^{2}+\mathcal{U}_{k,3}^{2}}{\left( P_{k,1} + \lambda_k \rho_k \zeta_1\right)^2}>\frac{1}{2} \h{20pt}\text{for any $k\in \mathbb{N}$ suitably large.}
\end{equation}
Here, $\mathcal{U}_{k}=\big(\mathcal{U}_{k,1},\mathcal{U}_{k,2},\mathcal{U}_{k,3}\big)^{\top}$. $P_{k, 1}$ is the first component of $P_k$. We will disprove \eqref{Contradiction equation in case 2} if (2) in \eqref{both and 2} holds. Preliminary lemmas are required before we show the proof. First, we establish the following relations between $\rho_{k}$, $\tau_{k}$ and the components of $z_{k}$.
\begin{armalem}\label{relation between r_{k} and e_{k}, case 2}The sequences $\{\rho_k\}$ and $\{\tau_k\}$ satisfy the following limits\h{0.3pt}$\mathrm{:}$
    $$\lim_{k \to \infty}  \tau_{k} = \lim_{k \to \infty}  \rho_{k} \h{0.3pt}\tau_{k}^{-2} = 0.$$ Denote by $z_{k,i}$ the $i$-th component of $z_{k}$, then
\begin{align*}
\frac{| z_{k,i} |}{\tau_{k}}\h{1.5pt}\lesssim\h{1.5pt} \frac{| P_k |}{\lambda_{k} \rho_{k}}, \h{20pt}\text{where $i=1,3$.}
\end{align*}
\end{armalem}
\begin{proof} First, we note that \begin{align}\label{set rel} D_{\lambda_k \rho_k}^+\left(P_k\right) \times \big\{ \phi : |\h{0.5pt} \phi \h{0.5pt}| < \pi \big\} \subset B^+_{ \left(1 + 2^{-4} \lambda_k \right) \h{0.5pt}| P_k |} \subset B^+_{ 2 \h{0.5pt}| P_k |} \subset B_{R_{\varepsilon_k}}^+.\end{align} It then turns out \begin{align*}
     \pi  \left| P_k \right| \int_{D_{\lambda_k \rho_k}^{+}(P_k)}  \left\vert D \tilde{u} \h{.8pt}\right\vert^{2}+\frac{4\tilde{u}_{1}^{2}+\tilde{u}_{3}^{2}}{\rho^2} \leqslant 2 \pi \int_{D_{\lambda_k \rho_k}^{+}(P_k)} \rho \left\vert D \tilde{u} \h{.8pt}\right\vert^{2}+\frac{4\tilde{u}_{1}^{2}+\tilde{u}_{3}^{2}}{\rho}   \leqslant \int_{B^+_{2 \h{0.5pt} \left| P_k \right|}} \left| \nabla \tilde{w} \right|^2.
\end{align*}This estimate induces $\tau_k^2 \leq \Theta\h{0.3pt}\big(\tilde{w}, 0, 2 \left| P_k \right|\big) \leqslant \varepsilon_k  \to 0$ as $k \to \infty$. Furthermore, (2) in \eqref{both and 2} yields
\begin{align*}
\rho_{k}^{\frac{1}{2}}\leqslant \Xi \h{0.3pt}\big(\tilde{u}, P_k, \theta_2\lambda_{k}\rho_{k}\big) \leqslant \tau_{k}^{2}.
\end{align*} Thus, $\rho_{k} \tau_{k}^{-2}\leqslant \tau_{k}^{2}\to 0$ as $k\to \infty$. The two limits in the lemma are obtained. \vspace{0.4pc}

By the definition of $z_k$, it holds, for $i = 1, 3$, that
\begin{align*}
| z_{k,i} | \h{1pt}\leqslant \h{1pt}\fint_{D_{1}^{+}} \big|\h{0.5pt} \mathcal{U}_{k,i} \h{0.5pt} \big| \h{1.5pt}\lesssim\h{1.5pt}  | P_k | \h{1pt} \left(\int_{D_{1}^{+}}\frac{\mathcal{U}_{k,i}^{2}}{\left(P_{k, 1} + \lambda_k \rho_k \zeta_1 \right)^{2}}\right)^{\frac{1}{2}}.
\end{align*}Through a change of variables and using the definition of $\tau_k$, we obtain \begin{align*}
| z_{k,i} | \h{1pt}\lesssim \h{1pt} \frac{| P_k |}{\lambda_k \rho_k} \big( \Xi \left(\tilde{u}, P_k, \lambda_k \rho_k \right) \big)^{\frac{1}{2}}.
\end{align*}The estimate of $z_{k, i}$ in the lemma then follows.
\end{proof}

\begin{armadef} We define a configuration space for the limit of $\big\{\mathcal{U}^{rs}_k\big\}$ as follows:
\begin{align*}
\mathcal{F}_{2}\coloneq \left\{\h{0.5pt}\mathcal{U}\in H^{1}\left(D_{\gamma_2}^{+};\mathbb{R}^{3}\right): \mathcal{U}=\mathcal{U}_{\infty}^{rs}\,\,\text{on}\,\,\partial^{+}D^+_{\gamma_2}\h{0.5pt}\right\}.
\end{align*}
\end{armadef}
Then, we claim that

\begin{armalem}\label{second H1 strong convergence prop}
As $k \to \infty$, $\left\{\mathcal{U}_{k}^{rs}\right\}$ converges to $\mathcal{U}^{rs}_{\infty}$ strongly in $H^{1}\big(D_{\gamma_2}^{+}\big)$. In addition, $\mathcal{U}^{rs}_{\infty}$ minimizes Dirichlet energy on $D_{\gamma_2}^{+}$ in the configuration space $\mathcal{F}_{2}$.
\end{armalem}
\begin{proof} We divide the proof into 3 steps. \vspace{0.3pc}

\noindent\textbf{Step 1. Comparison map.} Fix $k \in \mathbb{N}$, $R > 0$, $\mathcal{U} \in \mathcal{F}_{2}$ and define $$ \displaystyle
N_{k,R, \mathcal{U}}\coloneq
 z_{k}+R \tau_{k}\frac{\mathcal{U}}{\vert \mathcal{U} \vert \vee R}.$$
With $R$ and $\mathcal{U}$ fixed, it follows that \begin{align}\label{unif conv of N}
\left\| N_{k,R, \mathcal{U}} - z_k \right\|_{L^\infty(D_{\gamma_2}^+)} \leqslant R \tau_k \to 0, \h{20pt}\text{as $k \to \infty$.}
\end{align} Since $|z_k| \to 1$ as $k \to \infty$, the normalized vector field of $N_{k,R, \mathcal{U}}$, that is, $\widehat{N}_{k,R, \mathcal{U}} $, is well-defined for large $k$. The largeness of $k$ depends only on $R$ and is independent of $\mathcal{U}$. Using the above notation and a given $s\in (0,1)$, we define
\begin{equation*} \mathcal{U}_{k,s,R}(\zeta) :=
\begin{cases}
    N_{k,R, \mathcal{U}} \left( \dfrac{\zeta}{1-s} \right), &\h{15pt} \zeta \in D_{\left(1-s\right)\h{0.5pt}\gamma_2}^{+}; \\[4mm]
    \dfrac{\gamma_2 - |\zeta|}{s \gamma_2} \h{1pt} N_{k,R, \mathcal{U}^{rs}_{\infty}} \big( \h{0.8pt} \gamma_2 \widehat{\zeta}\h{1pt} \big) + \dfrac{ |\zeta|- \left(1-s\right) \gamma_2}{s \gamma_2} \h{1pt}\mathcal{U}_{k} \big(\h{0.8pt} \gamma_2 \widehat{\zeta} \h{1pt} \big), &\h{15pt} \zeta \in \overline{D^{+}_{\gamma_2}} \setminus D_{\left(1-s\right)\h{0.5pt}\gamma_2}^{+}.
\end{cases}
\end{equation*}
The normalized vector field of $\mathcal{U}_{k,s,R}$, that is, $\widehat{\mathcal{U}}_{k,s,R}$, is our comparison map. \vspace{0.4pc}

\noindent\textbf{Step 2. Validity of the comparison map.} Note that $\mathcal{U}_k^{rs} \to \mathcal{U}_\infty^{rs}$ uniformly on $\p^+ D_{\gamma_2}^+$. Thus, $$ \big\Vert\h{0.5pt} \mathcal{U}_k - z_k \h{0.5pt} \big\Vert_{L^{\infty}\big(\p^+ D_{\gamma_2}^+ \big)} \h{1.5pt}\leqslant \h{1.5pt}\tau_k \left( 1 + \big\|\h{0.5pt} \mathcal{U}_\infty^{rs} \h{0.5pt}\big\|_{L^{\infty}\big(\p^+ D_{\gamma_2}^+\big)}\right), \h{15pt}\text{provided that $k$ is suitably large.} $$ $\widehat{\mathcal{U}}_{k,s,R}$ is well-defined for large $k$ by this estimate and \eqref{unif conv of N}. \vspace{0.2pc}

On the other hand, for any $\zeta \in \p^+ D_{\gamma_2}^+$, it holds $\mathcal{U}_{k, s, R}(\zeta) = \mathcal{U}_k\big(\gamma_2 \widehat{\zeta}\h{1pt}\big) = \tilde{u}\h{0.3pt}\big(P_k + \lambda_k \rho_k \gamma_2 \widehat{\zeta}\h{1pt}\big)$. Moreover, we have $\mathcal{U}_{k, s,R}(\zeta) = N_{k,R, \mathcal{U}^{rs}_{\infty}} \big( \h{0.8pt} \gamma_2 \widehat{\zeta}\h{1pt} \big)$ for any $\zeta \in \p^+ D_{\left(1 - s\right)\gamma_2}^+$. The trace of $N_{k,R, \mathcal{U}} \big(  \zeta/(1-s) \big)$ on $\p^+ D_{\left(1-s\right)\gamma_2}^+$ also equals $N_{k,R, \mathcal{U}^{rs}_{\infty}} \big( \h{0.8pt} \gamma_2 \widehat{\zeta}\h{1pt} \big)$. Hence, the map $\mathcal{U}_{k, s, R}$ is $H^1$-regular on $D_{\gamma_2}^+$.\vspace{0.4pc}

\noindent\textbf{Step 3. Energy comparison.} Define $\tilde{v}$ as $L \big[\h{1pt} \widehat{\mathcal{U}}_{k,s,R}\left(\frac{\cdot - P_k}{\lambda_{k} \rho_{k}}\right) \big]$ on $T_k := D_{\lambda_k \rho_k \gamma_2}^+\left(P_k\right) \times \big\{ \phi : |\h{0.5pt}\phi\h{0.5pt}| < \pi \big\}$. Since $\tilde{w}$ is an energy minimizing map of  $\mathfrak{F}_{\nu}$, it holds  \begin{align}\label{key ine}
&\int_{T_k} g^{\alpha\beta}  \partial_{s_{\alpha}} \tilde{w} \cdot \partial_{s_{\beta}} \tilde{w} \h{1.5pt}\diff V_g \leqslant \mathfrak{F}_{\nu}[\tilde{w}, T_k] \leqslant \mathfrak{F}_{\nu}[\tilde{v},T_k] \nonumber\\[1.5mm] & \h{25pt}:= \int_{T_k} g^{\alpha\beta}  \partial_{s_{\alpha}} \tilde{v} \cdot \partial_{s_{\beta}} \tilde{v} +\sqrt{2}\mu \big(1-3S[\tilde{v}]\big)  \diff V_{g}+\nu\int_{\overline{T_k} \h{1pt}\cap \h{1pt} \{ z = 0 \}}\big\vert \tilde{v}-\tilde{w}^{\circ}\big\vert^{2}\diff S_{g}.
\end{align}

We now estimate the terms on the left and right-hand sides above.\vspace{0.2pc}

\noindent\textbf{Estimate 3.1.} First, it can be computed that \begin{align*}
    g^{\alpha\beta}  \partial_{s_{\alpha}} \tilde{w} \cdot \partial_{s_{\beta}} \tilde{w} &= \left|\h{0.5pt} \p_\rho \tilde{u} \h{0.5pt}\right|^2 \left( g^{11} \cos^2 \phi + g^{22} \sin^2 \phi + g^{12} \sin 2 \phi\right) + \left|\h{0.5pt}\p_z \tilde{u}\h{0.5pt}\right|^2 g^{33}  \\[1.5mm]
    & + 2 \left( g^{13} \cos \phi + g^{23} \sin \phi\right) \p_\rho \tilde{u} \cdot \p_z \tilde{u}   + \frac{4 \tilde{u}_1^2 + \tilde{u}_3^2}{\rho^2} \left( g^{11} \sin^2 \phi + g^{22} \cos^2 \phi - g^{12} \sin 2 \phi\right).
\end{align*}
Thus, \begin{align*}
    &\tau_k^{-2} \int_{T_k} g^{\alpha\beta}  \partial_{s_{\alpha}} \tilde{w} \cdot \partial_{s_{\beta}} \tilde{w} \h{1.5pt}\diff V_g \geqslant \tau_k^{-2} \int_{- \pi}^{\pi} \diff \phi \int_{D_{\lambda_k \rho_k \gamma_2}^+\left(P_k\right)} 2 \left( g^{13} \cos \phi + g^{23} \sin \phi\right) \p_\rho \tilde{u} \cdot \p_z \tilde{u} \h{1pt} \sqrt{\det( g _{\alpha\beta})} \h{1pt}\rho   \\[1.5mm]
    &\h{18pt}+  \tau_k^{-2} \int_{- \pi}^{\pi} \diff \phi \int_{D_{\lambda_k \rho_k \gamma_2}^+\left(P_k\right)} \left(\h{1pt}\left|\h{0.5pt} \p_\rho \tilde{u} \h{0.5pt}\right|^2 \left( g^{11} \cos^2 \phi + g^{22} \sin^2 \phi + g^{12} \sin 2 \phi\right) + \left|\h{0.5pt}\p_z \tilde{u}\h{0.5pt}\right|^2 g^{33} \right)\sqrt{\det( g _{\alpha\beta})} \h{1pt}\rho.
\end{align*}Now we change variables by letting $(\rho, s_3) = P_k + \lambda_k \rho_k \zeta$, where $\zeta \in D_{\gamma_2}^+$. If we still use $(g_{\alpha\beta})$ to represent the metric after the change of coordinates, then the last estimate can be reduced to \begin{align*}
    &\tau_k^{-2} \int_{T_k} g^{\alpha\beta}  \partial_{s_{\alpha}} \tilde{w} \cdot \partial_{s_{\beta}} \tilde{w} \h{1.5pt}\diff V_g \\[1.5mm]
    &\geqslant \int_{- \pi}^{\pi}  \int_{D_{\gamma_2}^+} 2 \left( g^{13} \cos \phi + g^{23} \sin \phi\right) \p_{1} \mathcal{U}^{rs}_k \cdot \p_{2} \mathcal{U}^{rs}_k \h{1pt} \sqrt{\det( g _{\alpha\beta})} \left( P_{k,1} + \lambda_k \rho_k \zeta_1 \right)   \\[1.5mm]
    &+  \int_{- \pi}^{\pi}  \int_{D_{\gamma_2}^+} \left(\h{1pt}\left|\h{0.5pt} \p_1 \mathcal{U}^{rs}_k \h{0.5pt}\right|^2   \left( g^{11} \cos^2 \phi + g^{22} \sin^2 \phi + g^{12} \sin 2 \phi\right) + \left|\h{0.5pt}\p_2 \mathcal{U}^{rs}_k\h{0.5pt}\right|^2  g^{33} \right)\sqrt{\det( g _{\alpha\beta})} \left( P_{k, 1} + \lambda_k \rho_k \zeta_1 \right).
\end{align*} Here, $\p_1$ and $\p_2$ are the partial derivatives of $\zeta_1$ and $\zeta_2$, respectively. Note that $\rho_k \in \big(0, 2^{-4} |P_k|\h{0.5pt}\big)$. Using the lower semi-continuity and (2) in Lemma \ref{properties of coordinates transformation}, we deduce that \begin{align*}
    \liminf_{k \to \infty} \h{1pt} \left| P_k \right|^{-1} \tau_k^{-2} \int_{T_k} g^{\alpha\beta}  \partial_{s_{\alpha}} \tilde{w} \cdot \partial_{s_{\beta}} \tilde{w} \h{1.5pt}\diff V_g \geqslant 2 \pi \int_{D_{\gamma_2}^+} \left| D \mathcal{U}_{\infty}^{rs} \right|^2.
\end{align*}We also use the fact that $\left(g_{\alpha\beta}\right)$ converges to the identity matrix uniformly on $D_{\gamma_2}^+$ as $k \to \infty$.\vspace{0.4pc}

\noindent\textbf{Estimate 3.2.} We estimate the two potential integrals on the right-hand side of \eqref{key ine}. By the unit-length condition, \begin{align*}
    \sqrt{2}\mu \int_{T_k} 1-3S[\tilde{v}]  \diff V_{g}+\nu\int_{\overline{T_k} \h{1pt}\cap \h{1pt} \{ z = 0 \}}\big\vert \tilde{v}-\tilde{w}^{\circ}\big\vert^{2}\diff S_{g} \h{1.5pt}\lesssim \h{1.5pt}  \left| P_k \right| \left(\lambda_k \rho_k \right)^2 + \nu \left| P_k \right| \lambda_k \rho_k .
\end{align*}It then turns out by Lemma \ref{relation between r_{k} and e_{k}, case 2} that \begin{align*}
    \sqrt{2}\mu \left| P_k \right|^{-1} \tau_k^{-2} \int_{T_k} 1-3S[\tilde{v}]  \diff V_{g}+\nu \left| P_k \right|^{-1} \tau_k^{-2} \int_{\overline{T_k} \h{1pt}\cap \h{1pt} \{ z = 0 \}}\big\vert \tilde{v}-\tilde{w}^{\circ}\big\vert^{2}\diff S_{g} \h{1.5pt}\lesssim\h{1.5pt}  \left(\frac{\lambda_k \rho_k}{\tau_k} \right)^2 + \frac{\nu  \lambda_k \rho_k}{\tau_k^2},
\end{align*} which converges to $0$ as $k \to \infty$.\vspace{0.4pc}

\noindent\textbf{Estimate 3.3.} We are left to estimate the first term on the right-hand side of \eqref{key ine}. \vspace{0.2pc}

Denote by $\overline{u}$ the map $\widehat{\mathcal{U}}_{k,s,R}\left(\frac{\cdot - P_k}{\lambda_{k} \rho_{k}}\right)$. Through change of variables, it holds  \begin{align*}
    &\int_{- \pi}^{\pi} \int_{D_{\lambda_k \rho_k \gamma_2}^+\left(P_k\right)}  \left\{\h{1pt}\left|\h{0.5pt} \p_\rho \overline{u} \h{0.5pt}\right|^2 \left( g^{11} \cos^2 \phi + g^{22} \sin^2 \phi + g^{12} \sin 2 \phi\right) + \left|\h{0.5pt}\p_z \overline{u}\h{0.5pt}\right|^2 g^{33} \right\} \h{1.5pt}\sqrt{\det( g _{\alpha\beta})} \h{1pt}\rho\h{1pt} + \\[1.5mm]
    & \int_{- \pi}^{\pi}\int_{D_{\lambda_k \rho_k \gamma_2}^+\left(P_k\right)} \left\{\h{1pt}2 \left( g^{13} \cos \phi + g^{23} \sin \phi\right) \p_\rho \overline{u} \cdot \p_z \overline{u} \right\}  \sqrt{\det( g _{\alpha\beta})} \h{1pt}\rho = \\[1.5mm]
    &\int_{- \pi}^{\pi} \int_{D_{\gamma_2}^+}   \left|\h{0.5pt} \p_1 \widehat{\mathcal{U}}_{k,s,R} \h{0.5pt}\right|^2 \left( g^{11} \cos^2 \phi + g^{22} \sin^2 \phi + g^{12} \sin 2 \phi\right)  \h{1.5pt}\sqrt{\det( g _{\alpha\beta})} \h{1pt}\left(P_{k,1} + \lambda_k \rho_k \zeta_1\right) + \\[1.5mm]
    &  \int_{- \pi}^{\pi}\int_{D_{\gamma_2}^+} \left\{\h{2pt} \left|\h{0.5pt}\p_2 \widehat{\mathcal{U}}_{k,s,R}\h{0.5pt}\right|^2 g^{33} + 2 \left( g^{13} \cos \phi + g^{23} \sin \phi\right) \p_1 \widehat{\mathcal{U}}_{k,s,R} \cdot \p_2 \widehat{\mathcal{U}}_{k,s,R} \right\}  \sqrt{\det( g _{\alpha\beta})} \h{1pt}\left( P_{k, 1} + \lambda_k \rho_k \zeta_1 \right).
\end{align*}Similar arguments for \eqref{limt on big} infer that
\begin{align*}
\limsup_{s \to 0} \h{1pt}\limsup_{ R \to \infty}\h{1pt} \limsup_{k \to \infty} \h{1pt} \tau_k^{-2}\int_{D_{\left(1 - s \right)\gamma_2}^{+}} \big|\h{0.5pt} D \widehat{\mathcal{U}}_{k, s, R} \h{0.8pt}\big|^2 \leqslant \int_{D_{\gamma_2}^+}  \big| D  \mathcal{U} \big|^2.
\end{align*}Using similar derivations for \eqref{es on ann}, we obtain \begin{align*}
  \limsup_{s \to 0} \h{1pt}\limsup_{ R \to \infty}\h{1pt} \limsup_{k \to \infty} \h{1pt}  \tau_k^{-2}\int_{D_{\gamma_2}^+ \setminus D_{\left(1 - s \right)\gamma_2}^{+}} \big|\h{0.5pt} D \widehat{\mathcal{U}}_{k, s, R} \h{0.8pt}\big|^2 = 0.
\end{align*}
 In light of (2) in Lemma \ref{properties of coordinates transformation} and the above two estimates, \begin{align*}
     &\limsup_{s \to 0} \h{1pt}\limsup_{ R \to \infty}\h{1pt} \limsup_{k \to \infty} \h{1pt} \left| P_k \right|^{-1} \tau_k^{-2} \int_{T_k}  g^{\alpha\beta}  \partial_{s_{\alpha}} \tilde{v} \cdot \partial_{s_{\beta}} \tilde{v} - \frac{4 \overline{u}_1^2 + \overline{u}_3^2}{\rho^2} \left( g^{11} \sin^2 \phi + g^{22} \cos^2 \phi - g^{12} \sin 2 \phi\right) \diff V_g
 \end{align*}is not greater than $\displaystyle 2 \pi \int_{D_{\gamma_2}^+}  \big| D  \mathcal{U} \big|^2$. \vspace{0.2pc}

On $D_{\lambda_k \rho_k \left(1 - s\right)\gamma_2}^+\left(P_k\right) \times \big\{ \phi : |\h{0.5pt}\phi\h{0.5pt}| < \pi \big\}$, it satisfies \begin{align*}
    |\h{0.5pt} \overline{u}_1 | + |\h{0.5pt} \overline{u}_3 | \h{1pt}\lesssim\h{1pt} | z_{k, 1} | + | z_{k, 3} | + R \tau_k,
\end{align*}which infers, together with Lemma \ref{relation between r_{k} and e_{k}, case 2}, that \begin{align}\label{pot}
    &\left| P_k \right|^{-1} \tau_k^{-2} \int_{D_{\lambda_k \rho_k (1 - s)\gamma_2}^+\left(P_k\right) \times \big\{ \phi : |\h{0.5pt}\phi\h{0.5pt}| < \pi \big\}} \frac{4 \overline{u}_1^2 + \overline{u}_3^2}{\rho^2} \left( g^{11} \sin^2 \phi + g^{22} \cos^2 \phi - g^{12} \sin 2 \phi\right) \diff V_g \\[2mm]
    &\lesssim \h{1pt} \left| P_k \right|^{-1} \left(\lambda_k \rho_k\right)^2  \frac{|z_{k, 1}|^2 + | z_{k, 3}|^2 + R^2 \tau_k^2}{\tau_k^2} \h{1pt}\lesssim \h{1pt} \left| P_k \right|^{-1} \left(\lambda_k \rho_k\right)^2 \left( R^2 + \left| P_k \right|^2 \left(\lambda_k \rho_k \right)^{-2}\right) \longrightarrow 0 \h{5pt}\text{as $k \to \infty$}. \nonumber
\end{align}For $\zeta \in D_{\gamma_2}^+ \setminus D_{\left(1 - s\right)\gamma_2}^+$, we have \begin{align*}
\mathcal{U}_{k,s,R}(\zeta) = z_k+ \dfrac{\gamma_2 - |\zeta|}{s \gamma_2} \h{1pt} \left( R \tau_{k}\frac{\mathcal{U}_\infty^{rs}}{\vert \mathcal{U}_\infty^{rs} \vert \vee R} \big( \h{0.8pt} \gamma_2 \widehat{\zeta}\h{1pt} \big)\right) + \dfrac{ |\zeta|- \left(1-s\right) \gamma_2}{s \gamma_2} \h{1pt}\tau_k\h{1pt} \mathcal{U}^{rs}_{k} \big(\h{0.8pt} \gamma_2 \widehat{\zeta} \h{1pt} \big).  \end{align*} As $\mathcal{U}_k^{rs} \to \mathcal{U}_\infty^{rs}$ uniformly on $\p^+ D_{\gamma_2}^+$, the first and third components of $\mathcal{U}_{k, s, R}$ can be estimated by \begin{align*}
    \left|\h{1pt} \left(\mathcal{U}_{k,s,R}\right)_1 \right| + \left|\h{1pt} \left(\mathcal{U}_{k,s,R}\right)_3 \right| \h{1pt}\lesssim\h{1pt} \left| z_{k, 1}\right| + \left| z_{k, 3} \right| + \left( R + \big\Vert \mathcal{U}_\infty^{rs} \big\Vert_{L^\infty\left(\p^+ D_{\gamma_2}^+\right)} \right) \tau_k
\end{align*} on $T_{k, 1-s, 1} := T_k \setminus D_{\lambda_k \rho_k \left(1 - s\right)\gamma_2}^+\left(P_k\right) \times \big\{ \phi : |\h{0.5pt}\phi\h{0.5pt}| < \pi \big\}$. By the same derivation for \eqref{pot}, it follows that\begin{align*}
    &\left| P_k \right|^{-1} \tau_k^{-2} \int_{T_{k, 1-s, 1}} \frac{4 \overline{u}_1^2 + \overline{u}_3^2}{\rho^2} \left( g^{11} \sin^2 \phi + g^{22} \cos^2 \phi - g^{12} \sin 2 \phi\right) \diff V_g  \longrightarrow 0 \h{5pt}\text{as $k \to \infty$}.
\end{align*}

In summary, we obtain \begin{align*}\limsup_{s \to 0} \h{1pt}\limsup_{ R \to \infty}\h{1pt} \limsup_{k \to \infty} \h{1pt} \left| P_k \right|^{-1} \tau_k^{-2} \int_{T_k}  g^{\alpha\beta}  \partial_{s_{\alpha}} \tilde{v} \cdot \partial_{s_{\beta}} \tilde{v}  \h{1.5pt} \diff V_g \leq 2 \pi \int_{D_{\gamma_2}^+}  \big| D  \mathcal{U} \big|^2. \end{align*}

Now, we apply Estimates 3.1-3.3 above to \eqref{key ine} and conclude \begin{align*}
    \int_{D_{\gamma_2}^+} \left| D \mathcal{U}_{\infty}^{rs} \right|^2 \leqslant \int_{D_{\gamma_2}^+}  \big| D  \mathcal{U} \big|^2.
\end{align*}The proof of this lemma is complete.  \end{proof}

\begin{proof}[\textit{Proof of (2) in Lemma \ref{Campanato Criteria}}] For $i=1,3$, by Lemma \ref{relation between r_{k} and e_{k}, case 2}, we have the following estimate: \begin{align*}\int_{D_{\theta_{2}}^{+}} \left(\frac{\lambda_{k}\rho_{k}}{\tau_{k}}\right)^{2}\h{1pt}\frac{4 \,\mathcal{U}_{k,1}^{2}+\mathcal{U}_{k,3}^{2}}{\left( P_{k,1} + \lambda_k \rho_k \zeta_1\right)^2} \h{1pt}\lesssim\h{1pt} \theta_{2}^{2} +  \left(\frac{\lambda_{k} \rho_{k}}{\left| P_k \right|} \right)^2 \int_{D_{\theta_2}^2} \big| \h{0.5pt}\mathcal{U}_k^{rs} \big|^2 \h{20pt}\text{for any $k\in \mathbb{N}$ suitably large.} \end{align*}Using this estimate, uniform boundedness of $L^2$-norm of $\mathcal{U}_k^{rs}$ over $D_1^+$, and Lemma \ref{second H1 strong convergence prop}, we then can take $k \to \infty$ in \eqref{Contradiction equation in case 2} and obtain $\displaystyle \theta_2^2 + \int_{D_{\theta_2 }^{+}} \left\vert D \mathcal{U}^{rs}_{\infty} \right\vert^{2}  \gtrsim \frac{1}{2}.$
In addition, it holds
$\big\Vert D \mathcal{U}_{\infty}^{rs}\big\Vert_{L^{\infty}(B_{3/8}^{+})}\lesssim \big\Vert D \mathcal{U}^{rs}_{\infty}\big\Vert_{L^{2}(D_{\gamma_2}^{+})}\lesssim 1$ by the mean value identity for harmonic functions. In light of these two estimates,  $\theta_2$ must be bounded from below by a positive universal constant, which yields a contradiction if $\theta_2$ is taken suitably small. The smallness of $\theta_2$ is universal.
\end{proof}

\section{Boojums and Their Local Structures}\label{booj} This section is devoted to the proof of Items (2)-(4) in Theorem~\ref{main theorem 1}. First, we prove (2) of Theorem \ref{main theorem 1} in Section \ref{(2) of main 1}. Then we show the local structures of director fields near boojums in Section \ref{ls of boo}.

\subsection{Values of \texorpdfstring{$w_\nu$}{TEXT} at the poles}\label{(2) of main 1}
Using the finite energy condition and Proposition \ref{main theorem in section 2}, we have $$w_{\nu}(\mathcal{N})= - e_3\quad \text{or }\quad w_{\nu}(\mathcal{N}) = e_3, \h{15pt}\text{where $e_3 = (0,0,1,0,0)^{\top}$.}$$ Indeed, if $\nu$ is suitably large, only the first case above may occur with the given boundary data $w^\circ$.
\begin{armalem}\label{determine the sign of the third component}
There exists $\nu_0 > 0$ such that $w_{\nu}(\mathcal{N})=w_{\nu}(\mathcal{S})= - e_3$ for all $\nu \geqslant \nu_0$.
\end{armalem}\vspace{-\baselineskip}
\begin{proof} It suffices to prove $\tilde{w}(0) = - e_3$. \vspace{0.2pc}

First, we give a refined monotonicity inequality. Recall $v_r$ introduced in the proof of Lemma \ref{almost monotonicity for type 1} and note that $\tilde{w}^\circ \equiv - e_3$ in a neighborhood of $0$. By further adjusting $\rho_0$ suitably small, it turns out from Lemma \ref{energy comparison property} that
\begin{align*}  \frac{\diff}{\diff r} \left(\frac{1+C_{7}\h{.5pt}r}{r}  \left(\int_{B_r^+} \vert\nabla \tilde{w}\vert^{2} +  \frac{\nu}{2}  \int_{D_r} \left| \tilde{w} - \tilde{w}^\circ \right|^2\right) \h{1pt}\right)   +C_{7} \h{1pt}r \geqslant
\frac{1}{r} \int_{\partial^{+}B^+_{r}} \vert \partial_{r} \tilde{w} \vert^{2}, \end{align*}for almost every $r \in \left(0,\min\big\{\h{0.5pt}\rho_{0}, (2C_{6})^{-1}\big\}\h{0.5pt}\right)$. Here, we also use \eqref{energy comparison 2}. Thus, \begin{align}\label{mono in} \frac{1+C_{7}\h{.5pt}r}{r}  \left(\int_{B_r^+} \vert\nabla \tilde{w}\vert^{2} +  \frac{\nu}{2}  \int_{D_r} \left| \tilde{w} - \tilde{w}^\circ \right|^2\right) + \frac{C_7}{2} r^2\end{align} does not decrease on $\left(0,\min\big\{\h{0.5pt}\rho_{0}, (2C_{6})^{-1}\big\}\h{0.5pt}\right)$.\vspace{0.2pc}

Fix an arbitrary $\rho \in \left(0,\min\big\{\h{0.5pt}\rho_{0}, (2C_{6})^{-1}\big\}\h{0.5pt}\right)$ and take $\nu$ large enough such that $\nu^{-1} < \rho$. Then, by the monotonicity of \eqref{mono in}, it follows that\begin{align}\label{app of mo}
    \nu \int_{B_{\nu^{-1}}^+} \vert\nabla \tilde{w}\vert^{2} +  \frac{\nu^2}{2}  \int_{D_{\nu^{-1}}} \left| \tilde{w} - \tilde{w}^\circ \right|^2  + \frac{C_7}{2} \nu^{-2} \leqslant \frac{1+C_{7}\h{.5pt}\rho}{\rho}  \left(\int_{B_{\rho}^+} \vert\nabla \tilde{w}\vert^{2} +  \frac{\nu}{2}  \int_{D_\rho} \left| \tilde{w} - \tilde{w}^\circ \right|^2\right) + \frac{C_7}{2} \rho^2.
\end{align}Recall that $w = w_\nu$. In addition, as $\nu \to \infty$ and up to a subsequence, $w_\nu$ converges strongly in $H^1_{\mathrm{loc}}\left(\Omega\right)$ to some $w_\infty$. $\nabla w$ converges strongly in $L^2\left(\Omega\right)$ to $\nabla w_\infty$. $\sqrt{\nu}\left( w - w^\circ\right)$ converges strongly in $L^2\left(\p B_1\right)$ to $0$. If we denote by $\tilde{w}_\infty$ the transformed $w_\infty$ under the change of variables $\psi$ given in Section \ref{cha of va}, then it holds by taking $\nu \to \infty$ in \eqref{app of mo} that \begin{align}\label{lim enr}
    \limsup_{\nu \to \infty} \left( \nu \int_{B_{\nu^{-1}}^+} \vert\nabla \tilde{w}\vert^{2} +  \frac{\nu^2}{2}  \int_{D_{\nu^{-1}}} \left| \tilde{w} - \tilde{w}^\circ \right|^2 \right) \leqslant \frac{1+C_{7}\h{.5pt}\rho}{\rho}  \int_{B_{\rho}^+} \vert\nabla \tilde{w}_\infty \vert^{2} + \frac{C_7}{2} \rho^2.
\end{align}Note that $\tilde{w}_\infty$ is a minimizing map of the energy: \begin{align*}
    \mathfrak{F}_{0}[\tilde{v},{B_{\rho_0}^{+}}] := \int_{B_{\rho_0}^{+}} g^{\alpha\beta}  \partial_{s_{\alpha}} \tilde{v} \cdot \partial_{s_{\beta}} \tilde{v} +\sqrt{2}\mu \big(1-3S[\tilde{v}]\big)  \diff V_{g},
\end{align*}subject to $\tilde{v} = - e_3$ on the flat boundary of $B_{\rho_0}^+$. Hence, $\tilde{w}_\infty$ is smooth up to the flat boundary of $B_{\rho_0}^+$. We now take $\rho \to 0$ in \eqref{lim enr} and obtain \begin{align}\label{lim nu str}
    \lim_{\nu \to \infty} \left( \nu \int_{B_{\nu^{-1}}^+} \vert\nabla \tilde{w}\vert^{2} +  \frac{\nu^2}{2}  \int_{D_{\nu^{-1}}} \left| \tilde{w} - \tilde{w}^\circ \right|^2 \right) = 0.
\end{align}

Rescale $\tilde{w}$ by defining $\overline{w}\left(\zeta\right) := \tilde{w}\left(\zeta / \nu\right)$. \eqref{lim nu str} can be rewritten as \begin{align*}
    \lim_{\nu \to \infty} \left( \int_{B_{1}^+} \vert\nabla \overline{w}\vert^{2} +  \frac{1}{2}  \int_{D_{1}} \left| \overline{w} + e_3 \right|^2 \right) = 0.
\end{align*}In addition, the monotonicity of \eqref{mono in} infers that \begin{align*}\frac{1+C_{7}\h{.5pt}r \nu^{-1}}{r}  \left(\int_{B_r^+} \vert\nabla \overline{w}\vert^{2} +  \frac{1}{2}  \int_{D_r} \left| \overline{w} + e_3 \right|^2\right) + \frac{C_7}{2} r^2 \nu^{-2}\end{align*} does not decrease on $\left(0,1\right)$ for any $\nu$ suitably large. The same arguments as in Sections \ref{sma en}-\ref{prf of 2.3} can be applied to show that $\overline{w}$ is uniformly bounded in $C^\alpha\big(\overline{B_{1/2}^+}\big)$ for large $\nu$ and some $\alpha \in (0,1)$. Arzel\`{a}-Ascoli theorem then implies the uniform convergence of $\overline{w}$ to $- e_3$ on $B_{1/2}^+$ as $\nu \to \infty$. In particular, $\displaystyle \lim_{\nu \to \infty} \overline{w}(0) = \lim_{\nu \to \infty} \tilde{w}(0) = - e_3$. The proof is complete.
\end{proof}

\subsection{Local structures of boojums}\label{ls of boo}
In this section, we analyze the director field in a neighborhood of $\mathcal{N}$. To determine the director field near $\mathcal{N}$, we need to compare the three eigenvalues given in \eqref{eigenvalues}. First, we investigate the strict positivity of the first component of $u$, where $w = L[u]$. Here, $w$ and $u$ also depend on $\nu$.
\begin{armalem}\label{sing of u1}
$u_{1} > 0$ in $\omega := \big(\big(0,\infty)\times \mathbb{R}\big)\setminus \overline{\mathbb{D}}$, where $\mathbb{D}$ is the open unit disk in the $(\rho,z)$-plane.
\end{armalem}
\vspace{-\baselineskip}
\begin{proof} Define
$$\overline{u}=(\overline{u}_{1},\overline{u}_{2},\overline{u}_{3})^{\top}\coloneq \left(\vert u_{1}\vert,u_{2},u_{3}\right)^{\top},\quad \textup{where } u = \left( u_{1},u_{2},u_{3}\right)^{\top}.$$
By the standard Sobolev theory, it holds
\begin{align*}
\left| D u \right|^2= \left| D \overline{u} \right|^2 \h{15pt}\text{a.e. in } \omega.
\end{align*}
In addition, for the polynomial $P$ defined in \eqref{Polynomial P}, it satisfies $1-3P(u)\geqslant 1-3P(\overline{u})$. Using the triangle inequality, we have $
\left|\overline{u}_{1}-(u^{\circ})_{1}\right|\leqslant \left|u_{1}-(u^{\circ})_{1}\right|$ on $\big\{(\rho, z) : \rho^2 + z^2 = 1\big\}$. Hence, if $u$ is a minimizer of $E_{\nu}$, then so is $\overline{u}$.\vspace{0.2pc}

Since $u$ is smooth off the $z$-axis, by Serrin's maximum principle, either $|u_1| > 0$ or $u_1 \equiv 0$ in $\omega$. If $u_1 \equiv 0$, then by the boundary condition in \eqref{EL of w}, it turns out $u_1^\circ \equiv 0$ on $\big\{(\rho, z) : \rho^2 + z^2 = 1\big\}$, which contradicts our non-trivial assumption for $u_1^\circ$. We must have $|u_1| > 0$ in $\omega$. Note that $u_1$ cannot have a zero point in $\omega$. If there is a point in $\omega$ on which $u_1$ is negative, then $u_1 < 0 $ throughout $\omega$. Since $u$ and $\overline{u}$ are all minimizers of $E_\nu$, it follows \begin{align*}
    \int_{\omega}  \sqrt{2} \mu \big(1-3P(u)\big) \h{0.5pt}  \rho  + \nu \int_{-1}^1 |u - u^{\circ}|^2 = \int_{\omega}  \sqrt{2} \mu \big(1-3P(\overline{u})\big) \h{0.5pt} \rho  + \nu \int_{-1}^1 | \overline{u} - u^{\circ}|^2.
\end{align*}Using the fact that $1-3P(u)\geqslant 1-3P(\overline{u})$ again induces $$\displaystyle \int_{-1}^1 \left| u_1 - \big(u^\circ\big)_1 \right|^2 \leqslant \int_{-1}^1 \left| - u_1 - \big(u^\circ\big)_1 \right|^2, \h{15pt}\text{ equivalently $\displaystyle \int_{-1}^1 u_1 \big(u^\circ\big)_1 \geqslant 0$.}$$ This is a contradiction to $u_1 < 0$ in $\omega$, as $\big(u^\circ\big)_1$ is assumed to be non-negative and non-trivial on $\big\{(\rho, z) : \rho^2 + z^2 = 1\big\}$. The above arguments show that $u_1 > 0$ in $\omega$. The proof is complete.
\end{proof}

Using the sign of $u_1$ in $\omega$, we claim that \begin{armalem}
    For suitably large $\nu$, there exists an $\epsilon = \epsilon_\nu > 0$ suitably small such that \begin{align}\label{ord rela}\lambda_3 > \lambda_1 > \lambda_2 \h{15pt}\text{in $\omega \h{0.5pt} \cap \h{0.5pt} \mathbb{D}_\epsilon$.}\end{align} Here, $\mathbb{D}_\epsilon := \Big\{ (\rho, z) : \rho^2 + (z - 1)^2 < \epsilon^2 \Big\}$. $\lambda_1$, $\lambda_2$ and $\lambda_3$ are given in \eqref{eigenvalues}.
\end{armalem}
\begin{proof}
    To simplify the notation, we define $$\displaystyle b=b(u)\coloneq \frac{\sqrt{3}}{2}u_{1}+\frac{1}{2}u_{2}.$$ A direct calculation yields $$\displaystyle \lambda_{3}-\lambda_{1}=  \frac{\sqrt{3}}{2} b+ \frac{1}{2} \sqrt{1-b^{2}}.$$ By Lemma \ref{sing of u1} and the unit length condition of $u$, it holds $u_1 > 0$ and $u_2 > -1$ on $\omega$. Therefore, $\displaystyle b > - \frac{1}{2}$ on $\omega$, which yields $\lambda_3 > \lambda_1$ on $\omega$. \vspace{0.2pc}

    We can also calculate $$\displaystyle \lambda_1 - \lambda_2 = - \frac{\sqrt{3}}{2} b+ \frac{1}{2} \sqrt{1-b^{2}}.$$ Since $u = (0, -1, 0)^\top$ at $\mathcal{N}$ when $\nu$ is large, it follows $b < 0$ on $\omega \h{0.5pt} \cap \h{0.5pt} \mathbb{D}_\epsilon$, provided that $\epsilon > 0$ is small. Thus, we obtain $\lambda_1 > \lambda_2$ on $ \omega  \h{0.5pt} \cap \h{0.5pt} \mathbb{D}_\epsilon$. The proof is complete.
\end{proof}

From the last lemma, $\lambda_3$ is the largest eigenvalue of the matrix \eqref{matrix} in the neighborhood $ \omega \cap \mathbb{D}_\epsilon$. Therefore, the director field $\mathbf{d}$ can be represented by the normalized $\mathbf{n}$-field in $\omega \cap \mathbb{D}_\epsilon$, where
\begin{align}\label{director field general formula}
\mathbf{n}\left(\rho,\varphi,z\right)\coloneq \frac{\sqrt{2}}{2} \left( 1 + \frac{u_1 - \sqrt{3}u_2}{\sqrt{(u_1 - \sqrt{3}u_2)^2 + 4u_3^2}} \right) \mathbf{e}_{\rho} + \frac{\sqrt{2}u_3}{\sqrt{(u_1 - \sqrt{3}u_2)^2 + 4u_3^2}}\h{1pt} \mathbf{e}_{z}.
\end{align}
$u$ is smooth in $\omega \cap \mathbb{D}_\epsilon$ and equals $(0,-1,0)^{\top}$ at $\mathcal{N}$. Hence, the director field $\mathbf{d}= \widehat{\mathbf{n}}$ is well defined and smooth in $\omega \cap \mathbb{D}_\epsilon$. In addition, it satisfies the following limit for any $z_0 \in [1, 1 + \epsilon)$
\begin{align*}
\lim_{\substack{(\rho, z) \h{1pt}\to \h{1pt} (0, z_0)\\[.6mm]\left(\rho, z\right) \h{1pt}\in \h{1pt} \omega \h{1pt}\cap\h{1pt} \mathbb{D}_\epsilon}}\mathbf{d}\left(\rho,\varphi,z\right) = \mathbf{e}_{\rho},
\end{align*}
which completes the proof of Item (4) in Theorem~\ref{main theorem 1}.

\section{Asymptotic Behaviors of \texorpdfstring{$w$}{TEXT} at Far Field}\label{far field arg}

In the previous sections, we show the smoothness of $w$ in $\omega$. The singularities of $w$, if they exist, must be located on the $z$-axis and away from the two poles. In this section, we analyze the asymptotic behavior of $w$ as $|x|\to\infty$. The main result is \begin{armaprop}\label{Asymptotic behavior at far-field}
 The density $e(w)$ decays according to the algebraic rate $|x|^{-1}$ as $x \to \infty$. In addition,
\begin{align}\label{asy at far}
\lim_{r \h{0.5pt} \to \h{0.5pt} \infty} \| w - e_3 \|_{L^\infty(\Omega \h{1pt}\setminus \h{1pt} B_r)} = 0
\end{align}
\end{armaprop}

To prove this proposition, we first show a Bochner-type inequality.

\begin{armalem}\label{Bochner-type inequality}\,
Suppose $w$ is a minimizer of the energy functional $F_{\nu}$ within the configuration space $\mathcal{F}_{\nu}$. There exists a constant $C_{17}>0$ such that the bulk energy density of $w$, that is $$e(w)\coloneq \vert \nabla w \vert^{2}+\sqrt{2}\mu\left(1-3S[w]\right),$$ satisfies the following Bochner-type inequality at every regular point of $w$:
\begin{align}\label{Bochner-type inequality equation}
-\Delta e(w)\leqslant 4 \vert \nabla w \vert^{4} + 6\sqrt{2}\mu \big(1 -3 S[w] \h{0.5pt}\big)\h{0.5pt}\vert \nabla w \vert^{2} +6\sqrt{2}\mu\left(\nabla w_{i} \cdot \nabla w_{j}\right) \left( \nabla^2_{w_i w_{j}}  S[w] - \delta_{ij} \right).
\end{align}
\end{armalem}
\vspace{-\baselineskip}
\begin{proof} Plugging $f=\partial_{\alpha}w_{i}$ into the identity
$\Delta (f^{2})=2\vert \nabla f \vert^{2}+2f\Delta f$ and summing all indices yields
\begin{align*}
\Delta \h{0.2pt}\vert \nabla w \vert^{2}=2 |\h{0.2pt} \nabla^2 w \h{0.2pt}|^2 +2 \left(\partial_{\alpha}w_{i}\right)\Delta \left(\partial_{\alpha}w_{i}\right).
\end{align*}
Therefore, we obtain, by the Euler-Lagrange equation in \eqref{EL of w}, that
\begin{align*}
-\Delta \h{0.2pt}\vert \nabla w \vert^{2} \leqslant - 2\left(\partial_{\alpha}w_{i}\right)\Delta \left(\partial_{\alpha}w_{i}\right) = 2\vert \nabla w \vert^{4} -9\sqrt{2}\mu\h{0.8pt}S[w] \h{0.8pt}\vert \nabla w \vert^{2} + 3\sqrt{2}\mu \h{0.8pt}\left(\nabla w_{i} \cdot \nabla w_{j} \right) \h{0.8pt}\nabla^2_{w_iw_{j}} S[w].
\end{align*}For the potential term, we have
\begin{equation*}\begin{aligned}
- \sqrt{2}\mu \h{0.5pt}\Delta \left(1-3S[w]\right)= 3\sqrt{2}\mu \,\partial_{\alpha}\left(\nabla_{w_{i}}S[w]\cdot \partial_{\alpha} w_{i}\right) = 3\sqrt{2}\mu \left(\left(\nabla w_{i} \cdot \nabla w_{j} \right) \h{0.8pt}\nabla^2_{w_iw_{j}} S[w] +  \nabla_{w_{i}}S[w]\cdot\Delta w_{i} \right).
\end{aligned}
\end{equation*}
Still using the Euler-Lagrange equation in \eqref{EL of w} infers\begin{align*}
3\sqrt{2}\mu \h{0.5pt} \nabla_{w_{i}}S[w]\cdot\Delta w_{i}=&2\left(-\Delta w_{i}-\vert \nabla w \vert^{2}w_{i}+\frac{9}{\sqrt{2}}\mu S[w]w_{i}\right)\Delta w_{i} \leqslant 2\vert \nabla w \vert^{4}-9\sqrt{2}\mu \h{0.5pt} S[w] \h{0.5pt}\vert \nabla w \vert^{2}.
\end{align*}
Combining the estimates above, we obtain
\begin{align*}
-\Delta e(w)\leqslant 4 \vert \nabla w \vert^{4}+6\sqrt{2}\mu\left( \left(\nabla w_{i} \cdot \nabla w_{j}\right) \nabla^2_{w_i w_{j}}  S[w]-3 S[w] \h{0.5pt}\vert \nabla w \vert^{2}\right).
\end{align*}The proof is complete.
\end{proof}

Since $w\in\mathcal{F}_{\nu}$, for any $ \varepsilon>0$, we can choose $\overline{R}_\varepsilon > 2$ sufficiently large such that $\displaystyle  \int_{\Omega\h{1pt}\setminus\h{1pt} B_{\overline{R}_{\varepsilon}}}e(w) <\varepsilon$. Fix a point $x \in \Omega \setminus B_{4 \overline{R}_\varepsilon}$, the function $\displaystyle
    \frac{1}{r}\int_{B_{r}(x)}e(w)$ does not decrease on $r \in (0, \overline{R}_\varepsilon)$. If $\varepsilon$ is chosen suitably small, then the standard interior $\varepsilon$-regularity result infers the H\"{o}lder continuity of $w$ on $\Omega \setminus B_{4 \overline{R}_\varepsilon}$. Applying similar argument as in \cite{borchers1980analyticity}, we conclude that $w$ is smooth in $\Omega\setminus B_{4\overline{R}_{\varepsilon}}$ with \begin{align}\label{unif bd of grad}\Vert \nabla w \Vert_{L^\infty\big(\Omega\h{1pt}\setminus\h{1pt} B_{4\overline{R}_{\varepsilon}}\big)} \leqslant  C_{\varepsilon}, \h{15pt}\text{for some constant $C_{\varepsilon}>0$.} \end{align} We now prove Proposition \ref{Asymptotic behavior at far-field}.

\begin{proof}[\textit{Proof of Proposition \ref{Asymptotic behavior at far-field}}]
The Bochner-type inequality \eqref{Bochner-type inequality equation} and the bound of $\nabla w$ in \eqref{unif bd of grad} imply that
\begin{align*}
-\Delta e(w)\leqslant C_\varepsilon \h{1pt}e(w)\quad\text{in}\h{2pt} \Omega\setminus B_{4 \overline{R}_{\varepsilon}}.
\end{align*}
Applying the local-boundness estimate (Theorem 4.1 in \cite{han2011elliptic}) and the monotonicity of $\displaystyle
    \frac{1}{r}\int_{B_{r}(x)}e(w)$, we obtain
\begin{align*}\sup_{B_{1/2}(x)}  e(w) \leqslant  C_{\varepsilon} \int_{B_1(x)} e(w) \leqslant C_\varepsilon | x |^{-1} \int_{B_{\frac{|x|}{2}}(x)} e(w),\h{20pt}\text{for all $x\in \Omega\setminus B_{8 \overline{R}_{\varepsilon}}$.}
\end{align*}The energy density $e(w)$ decays according to the algebraic rate $|x|^{-1}$ as $x \to \infty$.\vspace{0.2pc}

Using Morrey's inequality, we have $$\| w - e_3\|_{L^\infty(B_{1/2}(x))} \lesssim \| w - e_3\|_{L^6(B_{1/2}(x))} + \| \nabla w \|_{L^6(B_{1/2}(x))}.$$ When $|x| \geqslant 16 \overline{R}_\varepsilon$, it holds that $B_{1/2}(x) \subset B_{|x|\h{0.5pt}/\h{0.5pt}2}(x) \subset \Omega\setminus B_{8 \overline{R}_{\varepsilon}}$. Thus, \begin{align*}
    \| w - e_3\|_{L^\infty(B_{1/2}(x))} &\lesssim \| w - e_3\|_{L^6(\Omega \h{1pt}\setminus \h{1pt} B_{\overline{R}_\varepsilon})} + \| \nabla w \|^{\frac{2}{3}}_{L^\infty(B_{1/2}(x))} \left( \int_{B_{1/2}(x)} | \nabla w |^2 \right)^{\frac{1}{6}}\\[1mm]
    &\lesssim \| \nabla w \|_{L^2(\Omega \h{1pt}\setminus \h{1pt} B_{\overline{R}_\varepsilon})} + C_\varepsilon | x |^{- \frac{1}{3}} \left(  \int_{\Omega \h{1pt}\setminus \h{1pt}B_{\overline{R}_\varepsilon}} e(w)\right)^{\frac{1}{2}} \leqslant \varepsilon^{\frac{1}{2}} + C_\varepsilon \h{1pt}\varepsilon^{\frac{1}{2}} \h{1pt}|x|^{- \frac{1}{3}}.
\end{align*}It immediately implies \eqref{asy at far}. The proof is complete.
\end{proof}

\begin{armarem} In light of Proposition \ref{Asymptotic behavior at far-field} and Lemma \ref{determine the sign of the third component}, there must be an odd number of singularities on the $z$-axis above (also below) the colloid if $\nu$ is large enough. The local profiles of the singularities can be read from \cite{yu2020disclinations}.
\end{armarem}

\section{Saturn Ring Disclination on the Equatorial Plane}\label{saturn}

In this section, we prove Theorem~\ref{main theorem 3} under the hypothesis there. \vspace{0.2pc}

Note that $w_\nu = L[u_\nu]$ is the minimizer of the energy $F_\nu$ in the configuration space $\mathcal{F}_{\nu, \mathcal{R}}$. Using the arguments in \cite{yu2020disclinations}, $w_\nu$ is smooth in $\overline{\Omega}$, except at finitely many points on the $z$-axis $l_z$. Since $u^\circ_1 \geqslant 0$ on $\p \Omega$ and is not constant, $u_{\nu, 1}$, the first component of $u_\nu$, can be strictly positive on $\Omega \setminus l_z$, in light of Serrin's maximum principle. Note that $u_3^\circ \geqslant 0$ on $\p \Omega \cap \big\{ z \geqslant 0 \big\}$ and is not constant. Using the same argument, we can have $u_{\nu, 3} > 0$ on $\left(\Omega \setminus l_z\right) \cap \big\{ z > 0 \big\}$.  \vspace{0.2pc}

Now, we show that $w_\nu$ cannot be isotropic on $T$, where $$T := \big\{(\rho,0):\rho > 1\big\}.$$ When $\nu$ is large enough, it must be uniaxial at some point on $T$. In fact, under the $\mathcal{R}$-symmetric hypothesis in Theorem~\ref{main theorem 3}, $u_{\nu, 3} \equiv 0$ on $T$. By \eqref{eigenvalues}, it turns out that
$$\begin{aligned}\lambda_{1}=-\frac{1}{2}\left(u_{\nu,1}+\frac{1}{\sqrt{3}} u_{\nu,2}\right),\qquad \lambda_{2}=&\frac{1}{4}\left\{u_{\nu,1}+\frac{1}{\sqrt{3}} u_{\nu,2}-\left|\h{1pt}u_{\nu,1}-\sqrt{3} u_{\nu,2}\right| \h{1.5pt}\right\},\\[1.5mm] \lambda_{3}=&\frac{1}{4}\left\{u_{\nu,1}+\frac{1}{\sqrt{3}} u_{\nu,2}+\left|\h{1pt}u_{\nu,1}-\sqrt{3} u_{\nu,2}\right| \h{1.5pt}\right\}\qquad \text{on } T.\end{aligned}$$
The three eigenvalues cannot all be 0 at some point on $T$, since otherwise $u=0$ at this point, which contradicts the regularity of $u_{\nu}$ and the unit-length condition. Therefore, $u_{\nu}$ must be uniaxial or biaxial on $T$. If $u_\nu$ is uniaxial at some point on $T$, then only one of the following two scenarios may occur:
\begin{align*}
\text{(1). If $\lambda_{1}=\lambda_{2}$, then $u_{\nu}=\left(\frac{\sqrt{3}}{2}, -\frac{1}{2}, 0\right)$;} \h{20pt} \text{(2). If $\lambda_{2}=\lambda_{3}$, then $u_{\nu}=\left(\frac{\sqrt{3}}{2}, \frac{1}{2}, 0\right)$}.
\end{align*}Here, we also use the strict positivity of $u_{\nu, 1}$ on $T$. We claim Case (2) above must occur on $T$ if $\nu$ is large. \vspace{0.2pc}

To verify this, it is sufficient to show the following for some $\nu_0 > 0$: \begin{align}\label{suff con}
    u_{\nu, 1} - \sqrt{3} \h{0.5pt} u_{\nu, 2} > 0 \h{15pt}\text{at some $(\rho_\nu, 0) \in T$, for all $\nu \geqslant \nu_0$.}
\end{align}Indeed, we note that at far field, $u_\nu$ approximately equals $(0, 1, 0)$. Thus, \begin{align}\label{far field}u_{\nu, 1} - \sqrt{3}\h{0.5pt}u_{\nu, 2} < 0\h{15pt}\text{at $(\rho, 0) \in T$ for large $\rho$ and all $\nu > 0$}.\end{align} \eqref{suff con} and the smoothness of $u_\nu$ on $T$ infer $$u_{\nu, 1} - \sqrt{3}\h{0.5pt}u_{\nu, 2} = 0 \h{15pt}\text{at some $(r_\nu, 0)$ on $T$, for all $\nu \geq \nu_0$.}$$ Here, $r_\nu > \rho_\nu$. Equivalently, Case (2) occurs at $(r_\nu, 0)$. To prove \eqref{suff con}, we use a contradictory argument. Suppose that there exists a positive sequence $\nu_{k} \to  \infty$ as $k\to \infty$, such that $u_{\nu_k, 1} - \sqrt{3} \h{0.5pt} u_{\nu_k, 2} \leqslant 0$ on $T$. Since up to a subsequence, it holds that $w_{\nu_k}$ converges strongly in $H_{\mathrm{loc}}^1(\Omega)$ to some $w_*$, we then apply the compactness of the trace operator and obtain that $w_{\nu_k}$ converges to $w_*$ almost everywhere on $T$, up to a subsequence. Thus, \begin{align}\label{cond con}u_{*, 1} - \sqrt{3}\h{0.5pt}u_{*, 2} \leqslant 0 \h{15pt}\text{ on $T$, where $w_* = L[u_*]$.}\end{align} However, $w_*$ minimizes the energy \begin{align*}
 F_0[w] := \int_{\Omega} \big\vert \nabla w \big\vert^{2}+ \sqrt{2}\mu \left(\h{0.5pt}1-3 \h{1pt} S[w] \h{0.5pt}\right)
\end{align*} in the configuration space: \begin{align*}
\Big\{w\in H_{\text{loc}}^{1}(\Omega;\mathbb{S}^{4}) : &\text{ $w$ is $\mathcal{R}$-symmetric and axially symmetric},\\[1mm] &F_{0}[w]<+\infty,\hspace{5pt} w \,\,\text{satisfies}\, \eqref{far field f2}, \hspace{5pt} w = w^\circ\, \text{on $\p \Omega$}\h{1pt}\Big\}.
\end{align*} The standard Schoen-Unlenbeck argument implies the regularity of $w_*$ on $\overline{\Omega}$, except at finitely many singularities on $l_z$. Therefore, we obtain $u_{*, 1} - \sqrt{3}\h{0.5pt}u_{*, 2} > 0$ near $(1, 0)$, since $u_1^\circ - \sqrt{3}\h{0.5pt}u_2^\circ > 0$ at $(1, 0)$. This contradicts \eqref{cond con}. \eqref{suff con} follows. \vspace{0.2pc}

In the remainder, we investigate the Saturn ring structure of the disclination at $(r_\nu, 0)$. Due to \cite{morrey2009multiple}, $u_\nu$ is real analytic on $T$. $u_{\nu, 1} - \sqrt{3}\h{0.5pt}u_{\nu, 2}$ can only take the value $0$ finitely many times on $T$. Otherwise, $u_{\nu, 1} - \sqrt{3}\h{0.5pt}u_{\nu, 2} \equiv 0$ on $T$, which contradicts the far field condition of $u_\nu$. Thus, we can assume that
\[
u_{\nu, 1} - \sqrt{3}\h{0.5pt}u_{\nu, 2} > 0 \quad \text{on } \Big\{ (\rho, 0) : \rho \in (r_\nu - \epsilon, r_\nu) \Big\} \quad \text{and} \quad u_{\nu, 1} - \sqrt{3}\h{0.5pt}u_{\nu, 2} < 0 \quad \text{on } \Big\{ (\rho, 0) : \rho \in (r_\nu, r_\nu + \epsilon) \Big\},
\]for some small $\epsilon > 0$. Meanwhile, $u_{\nu, 3}$ cannot equal $0$ in $D_\epsilon(r_\nu, 0) \setminus T$. By \eqref{eigenvalues}, it turns out that
\begin{equation}\label{3>2}
\lambda_3 > \lambda_2 \qquad \text{on } D_\epsilon(r_\nu, 0) \setminus T.
\end{equation}We can also compute that $$\lambda_2 - \lambda_1 = \frac{1}{4} \left\{ \sqrt{3} \left( \sqrt{3}\h{0.5pt}u_{\nu, 1} + u_{\nu, 2} \right) - \sqrt{\left( u_{\nu, 1} - \sqrt{3} \h{0.5pt} u_{\nu, 2} \right)^2 + 4  u_{\nu, 3}^2} \right\}.$$ By the regularity of $u_\nu$ on $D_\epsilon(r_\nu, 0)$ and the fact that $u_\nu(r_\nu, 0) = \big(\frac{\sqrt{3}}{2}, \frac{1}{2}, 0 \big)$, it then turns out
\begin{equation*}
\lambda_2 > \lambda_1 \qquad \text{on } D_\epsilon(r_\nu, 0), \text{ provided } \epsilon \text{ is suitably small.}
\end{equation*} This inequality and \eqref{3>2} induce that
\begin{equation*}
\lambda_3 > \lambda_2 > \lambda_1 \qquad \text{on } D_\epsilon(r_\nu, 0) \setminus T, \text{ provided } \epsilon \text{ is suitably small.}
\end{equation*}
By the above arguments, on $D_\epsilon(r_\nu, 0) \setminus \big\{(r_\nu, 0)\big\}$, $u_\nu$ is biaxial with $\lambda_3$ being the largest eigenvalue.\vspace{0.2pc}

On $D_\epsilon(r_\nu, 0) \setminus T$, the director field, that is the normalized eigenvector corresponding to $\lambda_3$, is given by \eqref{n-field formula}. Since $u_{\nu, 3} \equiv 0$ and $u_{\nu, 1}-\sqrt{3}\h{0.5pt}u_{\nu,2} > 0$ on $\big\{(\rho, 0) : \rho \in (r_\nu - \epsilon, r_\nu)\big\}$, it holds $\mathbf{d}\left[u_{\nu}\right] \equiv \mathbf{e}_{\rho}$ on $\big\{(\rho, 0) : \rho \in (r_\nu - \epsilon, r_\nu)\big\}$. However, $\mathbf{d}\left[u_{\nu}\right]$ is not continuous on $\Big\{ (\rho, 0) : \rho \in (r_\nu, r_\nu + \epsilon) \Big\}$. Fix an $\epsilon' \in (0, \epsilon)$ and denote by $x_\nu'$ the point $(r_\nu + \epsilon', 0)$. When approaching $x_\nu'$ along the lower part of $\partial D_{\epsilon'}(x_0)$, the component $u_{\nu, 3}$ remains strictly negative. Then it follows that
\begin{align*}
& \frac{\sqrt{2} \h{0.5pt} u_{\nu, 3}}{\sqrt{\left( u_{\nu, 1} - \sqrt{3} \h{0.5pt} u_{\nu, 2} \right)^2 + 4  u_{\nu, 3}^2}} \left( 1 + \frac{u_{\nu, 1} - \sqrt{3} \h{0.5pt} u_{\nu, 2}}{\sqrt{\left( u_{\nu, 1} - \sqrt{3} \h{0.5pt} u_{\nu, 2} \right)^2 + 4 u_{\nu, 3}^2}} \right)^{-1/2} \\[1.5mm]
&\h{30pt}= -\frac{\sqrt{2}}{2} \left( \left( u_{\nu, 1} - \sqrt{3} \h{0.5pt} u_{\nu, 2} \right)^2 + 4  u_{\nu, 3}^2 \right)^{-1/4} \left( \sqrt{\left( u_{\nu, 1} - \sqrt{3} \h{0.5pt}u_{\nu, 2} \right)^2 + 4  u_{\nu, 3}^2} - \left( u_{\nu, 1} - \sqrt{3} \h{0.5pt} u_{\nu, 2} \right) \right)^{1/2},
\end{align*}
where the above equality is evaluated on the lower part of $\partial D_{\epsilon'}(x'_\nu)$. Since $u_{\nu, 1} - \sqrt{3} \h{0.5pt} u_{\nu, 2} < 0$ at $x'_\nu$, it turns out that
\[
\frac{\sqrt{2} \h{0.5pt} u_{\nu, 3}}{\sqrt{\left( u_{\nu, 1} - \sqrt{3} \h{0.5pt} u_{\nu, 2} \right)^2 + 4  u_{\nu, 3}^2}} \left( 1 + \frac{u_{\nu, 1} - \sqrt{3} \h{0.5pt} u_{\nu, 2}}{\sqrt{\left( u_{\nu, 1} - \sqrt{3} \h{0.5pt} u_{\nu, 2} \right)^2 + 4 u_{\nu, 3}^2}} \right)^{-1/2} \longrightarrow -1,
\]as $(\rho, z) \to x'_\nu$ along the lower part of $\partial D_{\epsilon'}(x'_\nu)$. Thus, $\mathbf{d}\left[u_{\nu}\right] \to - \mathbf{e}_z$, as $(\rho, z) \to x'_\nu$ along the lower part of $\partial D_{\epsilon'}(x'_\nu)$. Using the same derivations, we have $\mathbf{d}\left[u_{\nu}\right] \to  \mathbf{e}_z$, as $(\rho, z) \to x'_\nu$ along the upper part of $\partial D_{\epsilon'}(x'_\nu)$. By \eqref{n-field formula}, $\mathbf{d}\left[u_{\nu}\right]$ has positive coefficient in front of $\mathbf{e}_\rho$ when $(\rho, z) \in \partial D_{\epsilon'}(r_\nu, 0) \setminus \{x'_\nu\}$. Therefore, when we start from $x'_\nu$ and rotate counterclockwise along $\partial D_{\epsilon'}(r_\nu, 0)$ back to $x'_\nu$, the director field $\mathbf{d}\left[u_{\nu}\right]$ continuously varies from $\mathbf{e}_z$ to $-\mathbf{e}_z$. Meanwhile, $\mathbf{d}\left[u_{\nu}\right]$ remains on the right-half part of the $(\rho, z)$-plane. The total variation of the angle of $\mathbf{d}\left[u_{\nu}\right]$ is $- \pi$, which shows that the topological degree of the director field around
$(r_\nu, 0)$ is $- 1/2$. \vspace{0.2pc}

We are left to investigate the asymptotic behavior of $\mathbf{d}\left[u_{\nu}\right]$ near $(r_\nu, 0)$. Fixing an $\alpha \in (0, 2\pi)$, we define $x_{\nu, \alpha}(t) := \left(r_{\nu}, 0\right)+t\left(\cos \alpha, \sin \alpha\right)$. Applying L'H\^{o}pital's rule infers \begin{align*}
    &\lim_{t \to 0} \frac{u_{\nu, 1}(x_{\nu, \alpha}(t)) - \sqrt{3}\h{0.5pt}u_{\nu, 2}(x_{\nu, \alpha}(t))}{u_{\nu, 3}(x_{\nu, \alpha}(t))} \\[1.5mm]&\h{15pt}= \lim_{t \to 0} \frac{\p_\rho \left(u_{\nu, 1} - \sqrt{3}\h{0.5pt}u_{\nu, 2}\right) \big|_{x_{\nu, \alpha}(t)} \cos \alpha + \p_z \left(u_{\nu, 1} - \sqrt{3}\h{0.5pt}u_{\nu, 2}\right) \big|_{x_{\nu, \alpha}(t)} \sin \alpha}{\p_\rho u_{\nu, 3} \big|_{x_{\nu, \alpha}(t)} \cos \alpha + \p_z u_{\nu, 3} \big|_{x_{\nu, \alpha}(t)} \sin \alpha}.
\end{align*}Using the odd symmetry of $u_{\nu, 3}$
and the even symmetry of $u_{\nu, 1}$, $u_{\nu, 2}$ with respect to the $z$-variable, we obtain $$\p_\rho u_{\nu, 3} \h{1pt} \big|_{(r_\nu, 0)} = 0 \qquad\text{and}\qquad \left(\p_z u_{\nu, 1} - \sqrt{3}\h{0.5pt}\p_z u_{\nu, 2}\right)\big|_{(r_\nu, 0)} = 0.$$ By the Hopf Lemma (see Corrollary 2.9 in \cite{linhan}), it follows that \begin{align}\label{sign of par z}\p_z u_{\nu, 3} \big|_{(r_\nu, 0)} > 0. \end{align} If we define $$\varkappa :=\frac{\partial_\rho u_{\nu,1}-\sqrt{3} \h{0.5pt} \partial_\rho u_{\nu,2}}{\partial_z u_{\nu,3}} \h{1.5pt}\bigg|_{\left(r_{\nu}, 0\right)},$$
then it follows that \begin{align*}
    \lim_{t \to 0} \frac{u_{\nu, 1}(x_{\nu, \alpha}(t)) - \sqrt{3}\h{0.5pt}u_{\nu, 2}(x_{\nu, \alpha}(t))}{u_{\nu, 3}(x_{\nu, \alpha}(t))}  = \varkappa \h{1pt}\cot \alpha.
\end{align*}Item (3) in Theorem \ref{main theorem 3} holds by the last limit, the sign of $u_{\nu, 3}$ above and below $\big\{z = 0\big\}$, and the expression of $\mathbf{d}\left[u_{\nu}\right]$ in  \eqref{n-field formula}. In light of \eqref{sign of par z}, to verify that $\varkappa \leqslant 0$, it is equivalent to show that \begin{align*} \partial_\rho u_{\nu,1}-\sqrt{3} \h{0.5pt} \partial_\rho u_{\nu,2} \leqslant 0 \h{15pt}\text{at $(r_\nu, 0)$,}\end{align*} which is a consequence of $u_{\nu, 1} - \sqrt{3}\h{0.5pt}u_{\nu, 2} > 0$ on $(r_\nu - \epsilon, r_\nu)$ and $u_{\nu, 1} - \sqrt{3}\h{0.5pt}u_{\nu, 2} < 0$ on $(r_\nu , r_\nu + \epsilon)$.

\backmatter

\section*{Declarations}
\begin{itemize}
\item \noindent \textbf{Funding:} Y. Yu is supported by the RGC of HK, grants  Nos. 14303723, 14306622, and 14304821.\vspace{0.3pc}
\item \noindent \textbf{Conflict of interest:} The authors declare that they have no conflict of interest.
\end{itemize}

\bibliography{sn-bibliography}%

\end{document}